\documentclass[a4paper,12pt]{amsart}
\usepackage[english]{babel}
\usepackage[T1]{fontenc}

\usepackage{enumerate}
\usepackage{comment}
\usepackage{multirow}

\usepackage[a4paper,margin=2.5cm]{geometry}
\usepackage{amsfonts,amsmath,amsthm,amssymb,amsxtra,calligra,mathrsfs}
\usepackage{graphicx}
\usepackage{tikz}
\usepackage{tikz-cd} 
\usepackage{xcolor}
\usepackage{upgreek}
\usepackage{graphicx}
\usepackage{mathtools}
\usepackage{extarrows}
\usepackage{faktor}
\usepackage{tensor}
\usepackage{enumitem}
\usepackage{hyperref}
\hypersetup{
    colorlinks=true, 
    citecolor=blue,
    linkcolor=magenta,
	bookmarksnumbered=true,
}

\newcommand{\bZ}{\mathbb{Z}}
\newcommand{\bP}{\mathbb{P}}
\newcommand{\bC}{\mathbb{C}}

\newcommand{\bQ}{{\mathbb Q}}

\newcommand{\bH}{{\mathbb H}}

\newcommand{\cP}{{\mathcal P}}

\newcommand{\cO}{{\mathcal O}}

\newcommand{\cI}{{\mathcal I}}

\newcommand{\cG}{{\mathcal G}}

\newcommand{\sheafHom}{\mathscr{H}\text{\kern -3pt {\calligra\large om}}\,}

\def\D{\Delta}

\DeclareMathOperator{\IH}{\mathrm{IH}}

\DeclareMathOperator{\Spec}{Spec}

\DeclareMathOperator{\Supp}{Supp}

\DeclareMathOperator{\Ext}{Ext}

\DeclareMathOperator{\Pic}{\mathrm{Pic}}

\DeclareMathOperator{\CH}{\mathrm{CH}}
\DeclareMathOperator{\ch}{\mathrm{ch}}
\DeclareMathOperator{\codim}{\mathrm{codim}}

\DeclareMathOperator{\pic}{Pic}

\DeclareMathOperator{\gr}{Gr}

\newcommand\im{\text{\rm Im}}

\newcommand\lra{\longrightarrow}

\newcommand{\cF}{{\mathcal F}}

\newcommand\cE{{\mathcal{E}}}

\newcommand\cL{\mathcal{L}}

\newcommand\cC{{\mathcal{C}}}

\newcommand\cK{{\mathcal{K}}}

\newtheorem{defn-pro}{Definition-Proposition}
\newtheorem{defn-thm}{Definition-Theorem}
\newtheorem{thm}{Theorem}[section]
\newtheorem{assu}[thm]{Assumption}
\newtheorem{lem}[thm]{Lemma}
\newtheorem{cor}[thm]{Corollary}
\newtheorem{pro}[thm]{Proposition}

\newtheorem{defn}[thm]{Definition}

\newtheorem{question}[thm]{Question}
\newtheorem{exmp}[thm]{Example}
\newtheorem{conj}[thm]{Conjecture}
\newtheorem{rem}[thm]{Remark}

\theoremstyle{remark}

\title{Asymptotic Behaviors of Moduli of One-dimensional Sheaves on Surfaces}

\author{ Fei Si }
\address{Beijing International Center for Mathematical Research, Peking University, No. 5 Yiheyuan Road Haidian District, Beijing, P.R.China 100871}
\email{sifei@bicmr.pku.edu.cn}

\author{Feinuo Zhang}
\address{Shanghai Center for Mathematical Sciences, Fudan University, Jiangwan Campus, Shanghai, 200438, China}
\email{fnzhang21@m.fudan.edu.cn}
\date{}

\begin{document}

\begin{abstract}
  In this paper, we study the asymptotic behaviors of the Betti numbers and Picard numbers of the moduli space $M_{\beta,\chi}$ of one-dimensional sheaves supported in a curve class $\beta$ on $S$ with Euler characteristic $\chi$.
  We determine the intersection cohomology Betti numbers of $M_{\beta,\chi}$ when $S$ is a del Pezzo surface and $\beta$ is sufficiently positive. As an application, we formulate a $P=C$ conjecture regarding the refined BPS invariants for local del Pezzo surfaces.
\end{abstract}
\maketitle

\tableofcontents

\section{Introduction}
\subsection{Motivation and background}
Throughout this paper, we work over the field $\bC$ of complex numbers. We denote by $\bZ$, $\bZ_{\geq0}$, $\bZ_{>0}$, $\bQ$ respectively the ring of integers,
the sets of nonnegative, positive integers and
the field of rational numbers.
All schemes are assumed to be of finite type over $\bC$. By a variety, we mean an integral, separated scheme and by a surface, we mean a smooth, projective $2$-dimensional variety. A curve is a projective scheme of pure dimension $1$. All sheaves are assumed to be coherent.

Let $S$ be a surface with an ample line bundle $H$ (we refer to such a pair $(S,H)$ as a polarized surface). Given $r\in\bZ_{\geq0}$, a line bundle $L$ with first Chern class $c_1(L)=c_1\in H^2(S,\bZ)$ and a cohomology class $c_2\in H^4(S,\bZ)\cong\bZ$, denote by $M(r,L,c_2)$ ({\it resp.} $M(r,L,c_2)^s$) the (coarse) moduli space parametrizing polystable ({\it resp.} stable) sheaves $F$ (with respect to $H$) on $S$ with
rank $r(F)=r$, determinant $\det(F)\cong L$ and the second Chern class $c_2(F)=c_2$. 
If we fix the first Chern class $c_1(F)=c_1$ rather than the determinant, then the moduli spaces are denoted by $M(r,c_1,c_2)$ and $M(r,c_1,c_2)^s$, respectively. When $r\geq1$, the discriminant of $F$ is defined as
$$\D(F)=2rc_2(F)-(r-1)c_1(F)^2 .$$

The above moduli spaces have played a crucial role in many areas and their geometry has been appealing to algebraic geometers for a long time.
When $r=1$, $M(1,L,c_2)$ is isomorphic to the Hilbert scheme $S^{[c_2]}$ parametrizing zero-dimensional closed subschemes in $S$ of length $c_2$, which is a smooth projective variety of dimension $2c_2$. Starting from Donaldson's generic smoothness result \cite{Don90}, studies on $M(r,L,c_2)$ for $r\geq 2$ have shown that $M(r,L,c_2)$ has better geometric properties as the discriminant 
$\D:=\D(F)$ $(F\in M(r,L,c_2))$
becomes larger. For example, when $\D$ is large enough, Gieseker-Li \cite{GL94} \cite{GiesekerLi96} and O'Grady \cite{O'Grady96} showed that  $M(r,L,c_2)$ is irreducible, normal and has generically smooth irreducible components of the expected dimension.

Some asymptotic phenomena also occur in the computation of Betti numbers of $M(r,L,c_2)$. G\"{o}ttsche has computed the Betti numbers of $S^{[m]}$ for $m\in\bZ_{\geq0}$ (\cite{Got90}), which in particular implies that the $i$-th Betti number $b_i(M(1,L,c_2))$ stabilizes as $c_2$ becomes sufficiently large. The first two Betti numbers of $M(2,L,c_2)^s$ have been computed by Li in \cite{Li97} when $\D$ is sufficiently large, which are
$$b_1(M(2,L,c_2)^s)=b_1(S) \quad \text{and}\quad b_2(M(2,L,c_2)^s)=1+\binom{b_1(S)}{2}+b_2(S),$$
matching $\mathop{\lim}\limits_{m\to\infty} b_1(S^{[m]})$ and $\mathop{\lim}\limits_{m\to\infty} b_2(S^{[m]})$, respectively.
In \cite{CW22}, Coskun and Woolf proposed the following conjecture for $M(r,c_1,c_2)$, which is closely related to the Atiyah-Jones conjecture \cite{AJ78} \cite{Taubes84} in gauge theory (see \cite[\S 8]{CW22}).
\begin{conj}[{\cite[Conjecture 1.1]{CW22}}]
\label{CWconj}
Assume $r\geq 1$. Then for each $i\in\bZ_{\geq0}$, there is a constant $N(i)$ depending on $i$ such that for every integer $k\leq i$,
$$b_k(M(r,c_1,c_2))=b_{k,\infty}$$ 
when $\D\geq N(i)$ holds, where $b_{k,\infty}:=\mathop{\lim} \limits_{c_2\to\infty}b_k(M(1,c_1,c_2))$. 
\end{conj}  
Using the wall-crossing and blowup formulae, Coskun and Woolf verified Conjecture \ref{CWconj} for all rational surfaces.
It is very natural to ask the following question, which was also suggested in \cite{CW22}.
\begin{question}
\label{quest}
    Are there similar stabilization results as that in Conjecture \ref{CWconj} on the Betti numbers of $M(r,L,c_2)$ or $M(r,c_1,c_2)$ when $r=0$?
\end{question}

In contrast to the positive rank case, there are few general results when $r=0$. Indeed, the rank zero case is not a parallel generalization of the positive rank case since the techniques for torsion free sheaves are usually not applicable to torsion sheaves, which makes the above question even more interesting.
To distinguish the rank zero case, we use the following notation instead.
Given an effective divisor $\beta$ on $S$ and an integer $\chi$, denote by $M_{\beta,\chi}$ the moduli space parametrizing polystable sheaves $F$ (with respect to $H$) on $S$ with $r(F)=0$, $\det(F)\cong\cO_S(\beta)$ and Euler characteristic $\chi(F)=\chi$. Such sheaves are supported on curves in the complete linear system $|\beta|$, so they are called one-dimensional sheaves. An essential difficulty caused by these sheaves is that their supports can be singular, even non-reduced. 

A main motivation of this paper is to answer Question \ref{quest} affirmatively under suitable conditions. To do so, we need to figure out how to replace the condition "$\D\gg0$". We observe that neither "$\D|_{r=0}=\beta^2\gg0$" nor "$\beta\cdot H\gg0$" is a reasonable condition. For example, let $\sigma: S_1\rightarrow \bP^2$ be the blowup of a point on $\bP^2$ and let $e_1$ be the exceptional divisor. Then by \cite[Proposition 3.10]{CGKT20}, $b_2(M_{\sigma^\ast \cO_{\bP^2}(d),1})=b_2(M_{ \cO_{\bP^2}(d),1}) \rightarrow 2 $  while $b_2(M_{\sigma^\ast \cO_{\bP^2}(d)-e_1,1})\rightarrow 3$ as $d \rightarrow \infty$   by our calculation .      In this paper, we give a formulation of a positivity condition "$\beta\gg0$" which should be the right one and we take the first step towards understanding the asymptotic behaviors of $M_{\beta,\chi}$.  For del Pezzo surfaces, a quantitative  condition for "$\beta\gg0$"  is given in Definition \ref{Ai}.  We also give an example in which the asymptotic irreducibility fails, again implying that the $r=0$ case is different from the positive rank case in nature.

We propose the following conjecture on the stable cohomology of $M_{\beta,\chi}$, as an analogue of Conjecture \ref{CWconj} in the rank zero case.

\begin{conj}
\label{conj}
For any given $i\in\bZ_{\geq0}$, $\chi\in\bZ$,
every fixed ample divisor $\beta_0$ and $\beta=n \beta_0$ ($n \in\bZ_{>0}$), 
there exists a constant $N(i,\beta_0,\chi)$ depending on $i$, $\beta_0$ and $\chi$ such that for each integer $k\leq i$,
the $k$-th intersection cohomology Betti number of $M_{\beta,\chi}$ is
$$\dim \IH^k(M_{\beta,\chi})=\lim_{m \to\infty}b_k(S^{[m]})$$
when $n \geq N(i,\beta_0,\chi)$, where $\IH^\ast(-)$ is the intersection cohomology with $\bQ$-coefficients.
\end{conj}

Another motivation for our work comes from the important role played by one-dimensional sheaves in enumerative geometry. The total space $X=\mathrm{Tot}(K_S)$ (also referred to as a local surface) of the canonical divisor $K_S$ is a non-compact Calabi-Yau threefold for a regular surface $S$. In this case, the moduli space of one-dimensional sheaves is closely related to enumerative invariants of $X$ defined via Gromov-Witten/Donaldson-Thomas/Pandharipande-Thomas theory (see for example \cite{PT10} \cite{MT18}). 

Meanwhile, there is a proper Hilbert-Chow morphism from $M_{\beta,\chi}$ to $|\beta|$ obtained by taking the Fitting support of a sheaf, which induces a perverse filtration on the cohomology ring of $M_{\beta,\chi}$.
An analogue of the $P=W$ conjecture \cite{dCHM12} with a different origin, called the $P=C$ conjecture, was formulated in \cite{KPS} for the moduli space of one-dimensional sheaves on the projective plane $\bP^2$. The $P=C$ conjecture provides an explanation to an asymptotic product formula for the refined BPS invariants for the local $\bP^2$ via the asymptotic Betti numbers of the moduli space of one-dimensional sheaves, which stimulates our interest in extending it to a larger generality.

\subsection{Main results}
Our main results in this paper are the following.
\begin{thm}[{\it cf}. Theorem \ref{thm:betti}]
\label{mainbetti}
Let $S$ be a del Pezzo surface. For any given $i\in\bZ_{\geq0}$, $\chi\in\bZ$, every effective divisor $\beta$ satisfying $(\mathrm{A}_i)$ as defined in Definition \ref{Ai} and each integer $k\leq i$, we have
    $$\dim \IH^k(M_{\beta,\chi})=\lim_{m \to\infty}b_k(S^{[m]}).$$
    In particular, Conjecture \ref{conj} is true for $S$.
\end{thm}
 This result generalizes \cite[Theorem 1.7]{Yuan23} for $\bP^2$ to all del Pezzo surfaces (see Example \ref{exmp:special}). Besides, it allows us to compute some Betti numbers of $M_{\beta,\chi}$ when the arithmetic genus $p_a(\beta)$ of $\beta$ is large, extending the result in \cite[Theorem 4.13]{CGKT20}.

The strategy of our proof is to relate two natural projective bundle maps from different relative Hilbert schemes, one maps onto $S^{[m]}$ when $\beta$ is sufficiently positive and the other maps onto the moduli space parametrizing one-dimensional sheaves supported on integral curves. The full support theorem of Migliorini-Shende \cite{MS2013} helps us connect different relative Hilbert schemes while another full support theorem of Maulik-Shen \cite{MS23} is used to connect moduli spaces $M_{\beta,\chi}$ for different $\chi$.

We also obtain a partial result on the stabilization of the Picard number of the moduli on an arbitrary surface. For simplicity, we consider the case when $\chi=1$ and  write $M_\beta$ for $M_{\beta,1}$.

\begin{thm}[{\it cf}. Theorem \ref{picardnum}]
\label{mainpic}
Let $(S,H)$ be any polarized surface. 
Then for every divisor satisfying $(\mathrm{P})$ as defined in Definition \ref{P}, we have the following relation of Picard numbers 
\begin{equation}
\label{rho}
    \rho(M_{\beta})\geq\rho(S)+1.
\end{equation}
\end{thm}
The above result generalizes \cite[Proposition 9.1]{CW22} to the rank zero case.
We prove it by using the determinant line bundles and constructing testing curves in $M_\beta$. The construction of testing curves is more subtle than that in the positive rank case since Hecke modifications at points may change the stability of a one-dimensional sheaf. The assumption on $\beta$ is made to ensure the existence of certain nodal curves in $|\beta|$.


\begin{rem}
We expect the inequality (\ref{rho}) to be an equality when $\beta$ is sufficiently positive, because by Conjecture \ref{conj}, the stable Picard number should be $\rho(S^{[m]})=\rho(S)+1$ ($m\geq2$).
    By Markman \cite{markman07}, for all rational surfaces and all K3 surfaces, if $M_{\beta}$ is smooth, then $H^*(M_{\beta},\bQ)$ is generated by the Künneth components of the Chern classes of a universal sheaf. In this case, the number of generators gives an upper bound of $\rho(M_{\beta})$ by $\rho(S)+1$ and our expectation is true.
\end{rem}

As an application, we formulate the following $P=C$ conjecture for all del Pezzo surfaces, extending the original form in \cite{KPS}. We also provide evidence for it.

\begin{conj}[{{\it cf.} Conjecture \ref{P=C}}]
Let $S$ be a del Pezzo surface. Suppose that $\beta\cdot H$ and $\chi$ are coprime. 
With the notation as in Section \ref{sec:app}, we have the following $P=C$ identity 
  \begin{equation*}
      P_\bullet H^{\ast}(M_{\beta,\chi})=C_\bullet H^{\ast}(M_{\beta,\chi}).
  \end{equation*}
\end{conj}

\begin{rem}
    In the final stages of preparation of this paper, we learned of very recent work \cite{PS24} of Pi-Shen on the P=C conjecture and refined BPS invariants for local $\bP^2$. Their project and ours started independently, but our methods and results have grown to share considerable similarities. For this reason, \cite{PS24} and the present paper will be superseded by a joint paper of the four authors, extending results of both papers. We are posting this version of the present paper for record.
\end{rem}

The rest of this paper is organized as follows. In Section \ref{sec:pre}, we collect some notations, definitions and facts  which will be useful later. Section \ref{sec:stabbetti} is addressed to study the relative Hilbert schemes and the non-integral locus in a linear system so that we can prove Conjecture \ref{conj} for del Pezzo surfaces. In Section \ref{sec:pic}, we construct three types of testing curves and bound the Picard numbers. Section \ref{sec:app} is mainly for the formulation of the $P=C$ conjecture. At the end, we present in Section \ref{sec:other} an example in which the moduli space is not asymptotically irreducible.


\subsection*{Acknowledgements}
The second author is grateful to her advisor, Professor Jun Li for his help in preparing for this work, and to Professor Claire Voisin for useful discussions. Both authors thank Professor Zhiyuan Li for helpful comments.
The first author is supported by NSFC grant 12201011.

\section{Preliminaries}
\label{sec:pre}
In this section, we review some properties concerning moduli spaces of our interest.

\subsection{Stability, moduli spaces and generators of the cohomology ring}
Let $(S,H)$ be a polarized surface. For a nonzero sheaf $E$ on $S$, denote by $P_{H,E}\in\bQ[m]$ the Hilbert polynomial of $E$ with respect to $H$ ($m\in\bZ_{\geq0}$)
$$P_{H,E}(m):=\chi(E\otimes H^{\otimes m})=\sum_{i=0}^{\dim E}a_i(H,E)\frac{m^i}{i!},$$
whose degree $\dim E\in\bZ_{\geq0}$ is exactly the dimension of the support of $E$, called the dimension of $E$.
The reduced Hilbert polynomial $p_{H,E}$ is the following monic polynomial
$$p_{H,E}:=\frac{P_{H,E}}{a_{\dim E}(H,E)}.$$
We say $E$ is pure (of dimension $n$) if every nonzero subsheaf of $E$ has dimension $n$.

\begin{defn}[\cite{Simpson94}]
\label{def:stab}
   A pure sheaf $F$ on $S$ is called semistable (with respect to $H$) if for any nonzero proper subsheaf $G\subset F$, 
   \begin{equation}
   \label{ssineq}
       p_{H,G}(m) \le  p_{H,F}(m)
   \end{equation}
   for $m \gg 0$.
   The above sheaf $F$ is called stable if the inequality (\ref{ssineq}) is strict. A semistable sheaf is called polystable if it is the direct sum of stable sheaves.
\end{defn}

The moduli space of torsion free semistable sheaves was first constructed by Gieseker \cite{Gieseker77} and Maruyama \cite{Maruyama78}.
Later, Simpson \cite{Simpson94} proved that for any fixed Hilbert polynomial $P\in\bQ[m]$, there exists a projective coarse moduli space $M_H(P)$ whose closed points are in bijection with polystable sheaves on $S$ with Hilbert polynomial $P$.
The moduli spaces $M(r,L,c_2)$, $M(r,c_1,c_2)$ and $M_{\beta,\chi}$ in the introduction can be respectively obtained as subschemes of $M_H(P_j)$ $(j=1,2,3)$ for a suitable $P_j$ according to the fixed data. We denote by $M$ any one of the above three moduli spaces.

A sheaf $\cE$ on $M\times S$ is called universal if for every scheme $T$ and every flat family $\cF$ of sheaves in $M$ over $T$, there exists a unique morphism $\pi_T:T\to M$ and a line bundle $B$ on $M$ such that 
$$(\pi_T\times 1_S)^*\cE\cong \cF\otimes p^*B,$$
where $p:M\times S\to M$ is the projection and $1_S:S\to S$ is the identity.
Note that two universal sheaves differ by a line bundle pulled back from $M$.
By the Künneth formula, the $i$-th Chern class of $\cE$
$$c_i(\cE)\in H^{2i}(M\times S,\bQ)=\bigoplus_{j+k=2i}H^j(M,\bQ)\otimes_\bQ H^k(S,\bQ).$$
If we fix a basis $\{\nu_l:\nu_l\in H^{\deg\nu_l}(S,\bQ),\,l\in\Lambda\}$ ($\Lambda$ is an index set) for $H^*(S,\bQ)$ as a $\bQ$-linear space, then we can write
$$c_i(\cE)=\sum_{l\in\Lambda}\mu_l\otimes_\bQ\nu_l$$
for uniquely determined $\mu_l\in H^{2i-\deg \nu_l}(M,\bQ)$. These $\mu_l$ are called the Künneth components of $c_i(\cE)$. The sub-$\bQ$-algebra generated by $\{\mu_l\}_{l\in\Lambda}$ is independent of the choice of $\{\nu_l\}_{l\in\Lambda}$.

We state the following theorems on the cohomology rings of moduli spaces.
\begin{thm}[\cite{Beauville95}]
\label{beauville}
Suppose $M$ is smooth projective moduli space parametrizing stable sheaves and 
$$\Ext^2(E,F)=0$$
for every $E,F$ in $M$.
If there exists a universal sheaf $\cE$ on $M\times S$, then $H^{\ast}(M,\bQ)$ is generated by the K\"{u}nneth components of the Chern classes of $\cE$.
\end{thm}
 
\begin{thm}(\cite{markman02}, \cite[Theorem 1, Theorem 2]{markman07}) \label{markman}
Assume $S$ is a rational Poisson surface (i.e., a rational surface with a nonzero section of $\cO_S(-K_S)$) and $r\geq0$.
If $M(r,c_1,c_2)$ is smooth, then the singular cohomology ring $H^\ast(M(r,c_1,c_2),\bZ)$ with $\bZ$-coefficients is torsion free and the cycle class map \[\mathrm{cl}:\CH^\ast(M(r,c_1,c_2)) \rightarrow H^\ast(M(r,c_1,c_2),\bZ)\] from the Chow ring of $M(r,c_1,c_2)$ (with $\bZ$-coefficients) to the cohomology ring is an isomorphism. 
\end{thm}

\subsection{One-dimensional sheaves and determinant line bundles}
From now on, we focus on sheaves supported on curves. If $E$ is a sheaf of dimension $1$, then 
$$P_{H,E}(m)=(c_1(E)\cdot H)m+\chi(E).$$
In Definition \ref{def:stab}, the inequality (\ref{ssineq}) that characterizes semistability is equivalent to
\begin{equation}
    \label{ss1dim}
    \mu(G)\leq\mu(F)
\end{equation}
when $F$ is a pure one-dimensional sheaf, where the slope $\mu(F)$ is defined as
$$\mu(F):=\frac{\chi(F)}{c_1(F)\cdot H}.$$

For any pure one-dimensinoal sheaf $F$, any surjection $\psi:E_0\to F$ from a locally free sheaf $E_0$ can be completed to a short exact sequence
$$0\to E_1\stackrel{\phi}\to E_0\stackrel{\psi}\to F\to0,$$
where $E_1$ is locally free with $r(E_1)=r(E_0)$. The Fitting support $\Supp(F)$ of $F$ is defined (as a closed subscheme of $S$) by the vanishing locus of the determinant of $\phi$, which is independent of the locally free resolution of $F$.

For a given effective divisor $\beta$ and $\chi\in\bZ$, the moduli space $M_{\beta,\chi}$ admits a Hilbert-Chow morphism 
$$h_{\beta,\chi}:M_{\beta,\chi} \rightarrow |\beta|$$ 
which sends a polystable sheaf $F\in M_{\beta,\chi}$ to its Fitting support $\Supp(F)$. 



Let $U\subset|\beta|$ be the open subset parametrizing integral curves. 
Then there is another modular interpretation of $h_{\beta,\chi}^{-1}(U)$ as the compactified Jacobian $\overline{J}^d_\pi$ of degree $d$ relative to the universal family  \[\pi:\cC_U\rightarrow U,\] 
of integral curves in $|\beta|$, where $\cC_U \subset U \times S$ is a closed subscheme  and  
\begin{equation*}
    d=\chi+\frac{\beta(\beta+K_S)}{2}. 
\end{equation*}

\begin{pro}
\label{prop:smsupp}
Suppose a general curve in $|\beta|$ is smooth and connected. Then the open subscheme $h_{\beta,\chi}^{-1}(U_0)$ of $M_{\beta,\chi}$ is smooth, where $U_0\subset|\beta|$ is the open subset parametrizing smooth curves. If in addition $h^1(\cO_S(\beta))=h^2(\cO_S(\beta))=0$, then the dimension of $h_{\beta,\chi}^{-1}(U_0)$ is
     \begin{equation}
         \label{dimsmsupp}
          \dim h_{\beta,\chi}^{-1}(U_0)=\beta^2+\chi(\cO_S).
     \end{equation}   
\end{pro}
\begin{proof}
    Since semistable sheaves on smooth curves with Euler characteristic $\chi$ are invertible sheaves of degree $d$, $h^{-1}(U_0)$ is the relative Picard scheme $\Pic^d(\cC_{U_0}/U_0)$, where $\cC_{U_0}:=\pi^{-1}(U_0)$. By \cite[Corollary 5.14, Proposition 5.19]{Kleiman05}, $\Pic^d(\cC_{U_0}/U_0)$ is smooth over $U_0$ of relative dimension $p_a(\beta)$ (the arithmetic genus of $\beta$). Thus $h_{\beta,\chi}^{-1}(U_0)$ is smooth. 
    
    If $h^1(\cO_S(\beta))=h^2(\cO_S(\beta))=0$, then by the the adjunction formula 
    and the Riemann-Roch formula 
    \begin{equation*}
        \begin{aligned}
    \dim h_{\beta,\chi}^{-1}(U_0)&=\dim|\beta|+p_a(\beta)\\
    &=\left(\chi(\cO_S(\beta))-1\right)+\left(\frac{\beta(\beta+K_S)}{2}+1\right)\\
    &=\left(\frac{\beta(\beta-K_S)}{2}+\chi(\cO_S)\right)+\frac{\beta(\beta+K_S)}{2}\\
    &=\beta^2+\chi(\cO_S),
        \end{aligned}
    \end{equation*}
 which completes the proof.   
\end{proof}


Now we consider the case when $\beta\cdot H$ and $\chi$ are coprime. In this case, semistable sheaves in $M_{\beta,\chi}$ are all stable. 
By \cite[Corollary 4.6.6, p. 119]{HL10}, there is a universal sheaf $\cE$ on $M_{\beta,\chi}\times S$.

Let $K^0(S)$ be the Grothendieck group of locally free sheaves on $S$. 
It is a commutative ring with $1=[\cO_S]$ where the multiplication is given by 
$$[F_1]\cdot[F_2]:=F_1\otimes F_2$$ for locally free sheaves $F_1$ and $F_2$ on $S$.
There is a group homomorphism from $K^0(S)$ to the Picard group of $M_{\beta,\chi}$ (\cite[Definition 8.1.1, p. 214]{HL10})
\begin{equation*}
   \lambda_{\cE}:\  K^0(S) \rightarrow \pic(M_{\beta,\chi}),\quad  \lambda_\cE(\mathfrak{a}):=\det(p_{!}(\cE \otimes q^\ast \mathfrak{a} )),
\end{equation*}
where $p:M_{\beta,\chi}\times S\to M_{\beta,\chi}$ and $q:M_{\beta,\chi}\times S\to S$ are the projections.

Since two universal sheaves differ by a line bundle pulled back from $M_{\beta,\chi}$, the homomorphism $\lambda_\cE$ is independent of the choice of $\cE$ when restricted to $K_\beta^0(S)$ (\cite[Lemma 8.1.2 {\it iv)}, p. 214]{HL10}). Here $K_\beta^0(S)\subset K^0(S)$ is the orthogonal complement of the class of $F\in M_{\beta,\chi}$ with respect to the bilinear form $(\mathfrak{a},\mathfrak{b})\mapsto \chi(\mathfrak{a}\cdot \mathfrak{b})$ on $K^0(S)\times K^0(S)$. We denote the restriction by
\begin{equation}
    \label{lambda}
    \lambda:K_\beta^0(S)\to \pic(M_{\beta,\chi}).
\end{equation}

For a sheaf $E$ on a scheme $X$ and a sheaf $G$ on a scheme $Y$, we write $E\boxtimes G$ for the sheaf $p_X^*E\otimes p_Y^*G$ on $X\times Y$, where $p_X:X\times Y\to X$ and $p_Y:X\times Y\to Y$ are the projections.

\begin{pro}\label{firstchern}
Suppose $\beta$ is base-point-free.
Then the first Chern class of $\cE$ is
\[c_1(\cE)=\ch_1(\cE)=\tilde{h}^\ast\ c_1(\cO_{|\beta|}(1) \boxtimes \cO_S(\beta)) \ \in \ \CH^1(M_{\beta,\chi} \times S)  \] 
 where $\tilde{h}:=h_{\beta,\chi}\times 1_S: M_{\beta,\chi} \times S \rightarrow |\beta| \times S$.
\end{pro}
\begin{proof}
  The argument is similar to that of \cite[Lemma 2.1]{PS23}. Clearly, $\cE$ is a torsion sheaf (on $M_{\beta,\chi} \times S$) supported on $\tilde{h}^{-1}(\cC)$, the pullback of $\cC\subset|\beta| \times S$ along $\tilde{h}$, where $\cC$ is the universal family of curves in $|\beta|$. Note that since $\beta$ is base-point-free, $\cC$ is flat over $|\beta|$ and therefore $\tilde{h}^{-1}(\cC)$ is a divisor on $M_{\beta,\chi}\times S$ corresponding to the pullback of the line bundle
  \begin{equation}
  \label{univcurv}
      \cO_{|\beta| \times S}(\cC) \cong \cO_{|\beta|}(1) \boxtimes \cO_S(\beta),
  \end{equation}
  from which the result follows.
\end{proof}

\begin{cor}\label{cor:lambda0}
In the situation of Proposition \ref{firstchern},
    let $x\in S$ be a closed point and denote by $[\cO_x]\in K^0(S)$ the class of the structure sheaf of $x$. Then $[\cO_x]\in K_\beta^0(S)$ and
    $$\lambda([\cO_x])=h_{\beta,\chi}^*\cO_{|\beta|}(1).$$
\end{cor}
\begin{proof}
The assertion that $[\cO_x]\in K_\beta^0(S)$ follows immediately from a locally free resolution of $F\in M_{\beta,\chi}$.
    By \cite[Example 8.1.3 {\it i)}, p. 214]{HL10},
    $$\lambda([\cO_x])=p_*(\det(\cE)|_{M_{\beta,\chi}\times\{x\}}).$$
    It follows from Proposition \ref{firstchern} that
    $$p_*(\det(\cE|_{M_{\beta,\chi}\times\{x\}}))=p_*(p^*h_{\beta,\chi}^*\cO_{|\beta|}(1))=h_{\beta,\chi}^*\cO_{|\beta|}(1),$$
    which completes the proof.
\end{proof}

\subsection{Moduli on del Pezzo surfaces}
\label{subsec:delpezzo}
Let $S$ be a del Pezzo surface. Then $S$ is either $\bP^2$, $\bP^1\times\bP^1$ or the blowup $S_\delta$ of $\delta$ general points on $\bP^2$ $(\delta=1,2,\cdots,8)$. 
We will need the following properties.

\begin{lem}
\label{effdiv}
    If $D$ is a nef divisor on $S$, then
\begin{equation}
\label{vanishing}
    h^1(\cO_S(D))=h^2(\cO_S(D))=0.
\end{equation}
   In particular,
    $$h^0(\cO_S(D))=\chi(\cO_S(D))=\frac{D(D-K_S)}{2}+1.$$
\end{lem}
\begin{proof}
Since $D$ is nef and $-K_S$ is ample, $D-K_S$ is ample. Hence (\ref{vanishing}) follows from the Kodaira vanishing theorem \cite[Theorem 4.2.1, p. 248]{Lazars04}. Then $\chi(\cO_S(D))$ can be computed by the Riemann-Roch formula.
\end{proof}


In the following two results, $\beta\cdot H$ and $\chi$ are not necessarily coprime.

\begin{lem}
\label{smmoduli}
    The moduli space $M_{\beta,\chi}$ is irreducible and smooth at points corresponding to stable sheaves. In particular, $h_{\beta,\chi}^{-1}(U)$ is smooth, where $U\subset|\beta|$ is the open subset parametrizing integral curves.
\end{lem}
\begin{proof}
The irreducibility of $M_{\beta,\chi}$ is proved for toric del Pezzo surfaces in \cite[Theorem 2.3]{MS23} and generalized to all del Pezzo surfaces in \cite[Theorem 1.5]{Yuan23}.
    By \cite[Corollary 4.5.2, p. 113]{HL10}, it suffices to show $\Ext^2(F,F)=0$ for any stable sheaf $F\in M_{\beta,\chi}$. This is proved in \cite[Lemma 2.5]{MS23}. For any $E\in h_{\beta,\chi}^{-1}(U)$ and any nonzero proper subsheaf $G\subset E$, $c_1(G)=c_1(E)$ and the nonzero quotient $E/G$ is zero-dimensional. Therefore, $\chi(G)<\chi(E)$ and $E$ is stable by (\ref{ss1dim}).
\end{proof}

\begin{thm}[{\cite{MS23}},{\cite{Yuan23}}]
\label{thm:chiindep}
Suppose $\beta$ is base-point-free and that a general curve in $|\beta|$ is smooth and connected.
    Then there is an isomorphism
    \begin{equation}
    \label{chiindep}
        R(h_{\beta,\chi})_*\mathrm{IC}_{M_{\beta,\chi}}\cong\bigoplus_{l=0}^{2p_a(\beta)}\mathrm{IC}\left(\wedge^lR^1\pi_{0*}\bQ_{\cC_{U_0}}\right)[-l+p_a(\beta)]
    \end{equation}
    in the bounded derived category of mixed Hodge modules on $|\beta|$, where $\mathrm{IC}_{M_{\beta,\chi}}$ is the intersection complex of $M_{\beta,\chi}$ and  $\pi_0:\cC_{U_0}\to U_0$ is the restriction of $\pi:\cC_U\to U$ to the open subset $U_0\subset|\beta|$ parametrizing smooth curves.
\end{thm}

\begin{pro}
\label{oddvanish}
    When $\beta\cdot H$ are $\chi$ are coprime, the cohomology ring $H^*(M_{\beta,\chi},\bQ)$ is generated by the Künneth components of the Chern classes of a universal sheaf. In particular, for $i\in\bZ_{\geq0}$,
    $$b_{2i+1}(M_{\beta,\chi})=0.$$
\end{pro}
\begin{proof}
    By Theorem \ref{beauville} and Lemma \ref{smmoduli}, it remains to check $\Ext^2(E,F)=0$ for any $E,F\in M_{\beta,\chi}$, which is done in the proof of \cite[Lemma 2.5]{MS23}. The vanishing of odd Betti numbers follows easily from the definition of Künneth components and the fact that odd Betti numbers of $S$ vanish.
\end{proof}

\section{Stabilization of Betti numbers for moduli on del Pezzo surfaces}
\label{sec:stabbetti}
Throughout this section, $S$ is a del Pezzo surface.

\subsection{Relative Hilbert schemes of points}
For $n\in\bZ_{\geq0}$, denote by $\tilde{\pi}^{[n]}:\cC^{[n]}\to |\beta|$ the relative Hilbert scheme of $n$ points on the universal family $\cC$ of curves in $|\beta|$. For $C\in|\beta|$, the fiber of $\tilde{\pi}^{[n]}$ over $C$ is the Hilbert scheme $C^{[n]}$ parametrizing length $n$, zero-dimensional closed subschemes of $C$. Note that $\tilde{\pi}^{[0]}:|\beta|\to|\beta|$ is the identity and
$\tilde{\pi}:=\tilde{\pi}^{[1]}:\cC\to|\beta|$ is the natural projection.

The following notion of being $k$-very ample is a generalization of being base-point-free ($0$-very ample) and being very ample ($1$-very ample).

\begin{defn}
    We say a divisor $\beta$ on $S$ is $k$-very ample $(k\in\bZ_{\geq0})$ if given any zero-dimensional closed subscheme $Z\subset S$ of length $k+1$, the restriction map 
    $$r_Z:H^0(\cO_S(\beta))\to H^0(\cO_S(\beta)|_Z)$$ is surjective.
\end{defn}

Now we can formulate the condition $(\mathrm{A}_i)$ to make "$\beta\gg0$" rigorous.

\begin{defn}
\label{Ai}
For $i\in\bZ_{\geq0}$,
    we say a divisor $\beta$ on $S$ satisfies $(\mathrm{A}_i)$ if 
\begin{enumerate}
\item $\beta$ is $i$-very ample and
$$2\dim|\beta|\geq 3i+2;$$
\item a general curve in $|\beta|$ is smooth and connected;
\item the codimension $\operatorname{codim}(|\beta|\setminus U,|\beta|)$ of the locus of non-integral curves in $|\beta|$ satisfies 
$$2\operatorname{codim}(|\beta|\setminus U,|\beta|)>i+1.$$
\end{enumerate}
\end{defn}



Since all curves in $|\beta|$ lie in $S$, there is a natural morphism $\sigma_i:\cC^{[i+1]}\to S^{[i+1]}$ sending a closed subscheme $Z\subset C$ $(C\in|\beta|)$ of length $i+1$ to $Z\subset S$. The morphism $\sigma_i$ allows us to regard $\cC^{[i+1]}$ as a scheme over $S^{[i+1]}$. The following property of $\sigma_i$ also appears in \cite[Proposition 3.16]{CGKT20}.

\begin{pro}
\label{iveryamp}
    Suppose $\beta$ is $i$-very ample. Then $\cC^{[i+1]}$ is a projective bundle over $S^{[i+1]}$.
\end{pro}
\begin{proof}
Let $\cI$ be the universal ideal sheaf on $S^{[i+1]}\times S$. Denote by $p:S^{[i+1]}\times S\to S^{[i+1]}$ and $q:S^{[i+1]}\times S\to S$ the projections. For every zero-dimensional closed subscheme $Z\subset S$ of length $i+1$, we have a short exact sequence
$$0\lra H^0(\cO_S(\beta)\otimes\cI_Z)\lra H^0(\cO_S(\beta))\stackrel{r_Z}\lra H^0(\cO(\beta)|_Z)\lra 0$$
since $\beta$ is $i$-very ample.
Thus $h^0(\cO_S(\beta)\otimes\cI_Z)=h^0(\cO_S(\beta))-(i+1)$,
which is independent of the choice of $Z$.
By Grauert's theorem (\cite[Chapter III, Corollary 12.9, p. 288]{Hartshorne}), 
$$\cG:=p_*(\cI\otimes q^*\cO_S(\beta))$$ is a vector bundle on $S^{[i+1]}$ of rank $h^0(\cO_S(\beta))-(i+1)$. Then it is easy to show that 
$\cC^{[i+1]}\cong  \bP(\cG^\vee)$.
\end{proof}

In what follows, we assume that $\beta$ satisfies $(\mathrm{A}_0)$ as defined in Definition \ref{Ai}.
Let $\pi^{[n]}$ ({\it resp.} $\pi^{[n]}_0$) be the restriction of $\tilde{\pi}^{[n]}$ to the open subscheme $U$ ({\it resp.} $U_0$) of $|\beta|$ parametrizing integral ({\it resp.} smooth) curves and let 

\begin{equation}
\label{notation:cun}
    \cC_U^{[n]}:=(\tilde{\pi}^{[n]})^{-1}(U),\quad\cC_{U_0}^{[n]}:=(\tilde{\pi}^{[n]})^{-1}(U_0).
\end{equation}
Note that
$\pi^{[1]}=\pi:\cC_U\to U$ and $\pi_0^{[1]}=\pi_0:\cC_{U_0}\to U_0$ are the natural projections.

Migliorini and Shende proved the following full support theorem for relative Hilbert schemes.
\begin{thm}[{\cite[Theorem 1]{MS2013}}]
\label{thm:support}
If $\cC_U^{[n]}$ is smooth, then
 \begin{equation}
 \label{eq:support}
     R \pi^{[n]}_\ast \bC[n+\dim |\beta|] =\mathop{\bigoplus} \limits_{l=-n}^n \mathrm{IC}(R^{n+l}\pi_{0\ast}^{[n]}  \bC)[-l].
 \end{equation}
 where $\mathrm{IC}(R^{n+l}\pi_{0 \ast}^{[n]}\bC)$ is the intermediate extension of the local system $R^{n+l}\pi_{0 \ast}^{[n]}  \bC$ on $U_0$.
\end{thm}

The family $\pi_0:\cC_{U_0}\to U_0$ of smooth curves admits a multisection 
\begin{equation}
\label{multisec}
    D_0\subset\cC_{U_0}
\end{equation}
such that $D_0$ is finite and flat over $U_0$ of degree $d_0\in\bZ_{>0}$ (see for example \cite[\S 5.2, p. 42]{MSY23}). Then the operation of adding a $0$-dimensional closed subscheme of length $d_0$ given by $D_0$ yields a morphism (a similar construction is used in \cite[\S 5.1]{Rennemo18})
$$\iota:\cC_{U_0}^{[n]}\to\cC_{U_0}^{[n+d_0]}.$$
There is a natural map 
$$\tilde{\iota}:\bC_{\cC_{U_0}^{[n+d_0]}}\to \iota_*\bC_{\cC_{U_0}^{[n]}}$$
and it induces 
\begin{equation*}
    \tilde{r}:R^{i}\pi_{0*}^{[n+d_0]} \bC\to R^{i}\pi_{0*}^{[n]} \bC.
\end{equation*}
On the stalks over each $b\in U_0$, the map $\tilde{r}_b$ induced by $\tilde{r}$ is the restriction map induced by $\iota_b:=\iota|_{\cC_b^{[n]}}$
$$\iota_b^\ast=\tilde{r}_b:H^i(\cC_b^{[n+d_0]},\bC)\to H^i(\cC_b^{[n]},\bC),$$
where $\cC_b^{[k]}:=(\tilde{\pi}^{[k]})^{-1}(b)$ ($k\in\bZ_{\geq0}$).

\begin{lem}
\label{lem:isom}
   If $i\leq n-1$, then the map $\tilde{r}$ is an isomorphism.
\end{lem}
\begin{proof}
It suffices to show $\tilde{r}_b$ is an isomorphism for each $b\in U_0$. Note that the restriction $\iota_b:\cC_b^{[n]}\to\cC_b^{[n+d_0]}$ of $\iota$ is the composition of
$$\iota_b^j:\cC_b^{[n+j-1]}\to\cC_b^{[n+j]}\quad(j=1,\cdots,d_0),$$
where $\iota_b^j$ is the morphism defined by adding a point $x_j\in\cC_b$ and $\sum_{j=1}^{d_0}x_j$ is the $0$-dimensional scheme $(D_0)_b:=D_0\cap\pi_0^{-1}(b)$.
By the proof of \cite[Chapter \uppercase\expandafter{\romannumeral7}, Proposition 2.2, p. 310]{ACGH85}, $\iota_b^j$ realizes $\cC_b^{[n+j-1]}$ as an ample divisor on $\cC_b^{[n+j]}$. By Lefschetz hyperplane theorem \cite[Theorem 3.1.17, p. 192]{Lazars04}, the restriction map
    $$\tilde{r}_b^{\,j}:H^i(\cC_b^{[n+j]},\bC)\to H^i(\cC_b^{[n+j-1]},\bC)$$ is an isomorphism since $i\leq n-1\leq n+j-2$. Then $\tilde{r}_b$ is an isomorphism since it is the composition of $\tilde{r}_b^{\,j}$.
\end{proof}

Now we can prove the following property relating the Betti numbers of the relative Hilbert schemes $\cC_U^{[n]}$ of integral curves for different $n$.

\begin{pro}
\label{n+d0}
    If $i\leq n-1$, then there is an isomorphism of $\bC$-vector spaces
    $$H^i(\cC_U^{[n]},\bC)\cong H^i(\cC_U^{[n+d_0]},\bC).$$
\end{pro}
\begin{proof}
    By \cite[Proposition 14]{She12} and Lemma \ref{smmoduli}, the relative Hilbert schemes $\cC_U^{[n]}$ are smooth for all $n\in\bZ_{\geq0}$. Thus we can use Theorem \ref{thm:support}. Applying the hypercohomology functor $\bH^{i-n-\dim|\beta|}(-)$ to both sides of (\ref{eq:support}), we obtain
\begin{equation*}
\label{eq:[n]}
\begin{aligned}
     H^i(\cC_U^{[n]},\bC)&= \bigoplus_{l=-n}^n \bH^{i-n-\dim|\beta|-l}\left(\mathrm{IC}(R^{n+l}\pi_{0\ast}^{[n]}\bC)\right)   \\
     &=\bigoplus_{l=-n}^{i-n}\bH^{i-n-\dim|\beta|-l}\left(\mathrm{IC}(R^{n+l}\pi_{0\ast}^{[n]}\bC)\right),
\end{aligned}
\end{equation*}
where the second equality holds since $\bH^{i-n-\dim|\beta|-l}\left(\mathrm{IC}(R^{n+l}\pi_{0\ast}^{[n]}\bC)\right)=0$
for $l>i-n$ by the vanishing result for perverse sheaves in \cite[Proposition 8.6.11, p. 148]{Maxim19}.
Similarly, 
\begin{equation*}
\label{eq:[n+d0]}
\begin{aligned}
    H^i(\cC_U^{[n+d_0]},\bC)&=\bigoplus_{m=-(n+d_0)}^{i-(n+d_0)}\bH^{i-(n+d_0)-\dim|\beta|-m}\left(\mathrm{IC}(R^{n+d_0+m}\pi_{0\ast}^{[n+d_0]}\bC)\right)\\
    &=\bigoplus_{l=-n}^{i-n}\bH^{i-n-\dim|\beta|-l}\left(\mathrm{IC}(R^{n+l}\pi_{0\ast}^{[n+d_0]}\bC)\right).
\end{aligned}
\end{equation*}
By Lemma \ref{lem:isom}, $R^{n+l}\pi_{0\ast}^{[n+d_0]}\bC\cong R^{n+l}\pi_{0\ast}^{[n]}\bC$ for $l=-n,\cdots,i-n$ and therefore 
each direct summand of $H^i(\cC_U^{[n]},\bC)$ is isomorphic to a direct summand of $H^i(\cC_U^{[n+d_0]},\bC)$. Hence the result follows.
\end{proof}

\subsection{Non-integral locus in a linear system}

Now we estimate the codimension of the non-integral locus in $|\beta|$ when $\beta$ satisfies $(\mathrm{A}_0)$ as defined in Definition \ref{Ai}. Given a non-integral curve $C\in|\beta|$, suppose it has irreducible components $C_1,\cdots, C_{m}$ ($m\in\bZ_{>0}$) with multiplicities $n_1,\cdots,n_{m}\in\bZ_{>0}$, respectively. 

If $m=1$, then $n_1\geq2$ and $C$ is an irreducible, non-reduced curve. By Lemma \ref{effdiv}, if $C_1$ is nef, then so is $\beta$ and
\begin{equation*}
\begin{aligned}
    \dim|\beta|-\dim|C_1|&=\frac{\beta(\beta-K_S)}{2}-\frac{n_1^{-1}\beta(n_1^{-1}\beta-K_S)}{2}\\
    &=(1-n_1^{-1})\frac{\beta(\beta-K_S)}{2}+\frac{n_1-1}{2n_1^2}\beta^2\geq\frac{\beta(\beta-K_S)}{4}.
\end{aligned}   
\end{equation*}
Let $Z_1$ be the locus of irreducible, non-reduced curves in $|\beta|$. Then it follows that
\begin{equation}
    \label{nonint1}
    \codim(Z_1,|\beta|)\geq \frac{\beta(\beta-K_S)}{4}.
\end{equation}

If $m\geq2$, then $C$ is reducible. If none of its irreducible components is a $(-1)$-curve, then $C_1,\cdots, C_{m}$ are all nef divisors. The locus $Z_2$ of such $C$ is contained in the image of
$|\beta_1|\times |\beta_2|$ in $|\beta|$, where $\beta_1$ and $\beta_2$ are nonzero, nef, effective divisors in $S$ with $\beta_1+\beta_2=\beta$. There are finitely many choices for $(\beta_1,\beta_2)$. By Lemma \ref{effdiv},
\begin{equation*}
    \dim|\beta|=\frac{(\beta_1+\beta_2)(\beta_1+\beta_2-K_S)}{2}=\dim|\beta_1|+\dim|\beta_2|+\beta_1\cdot\beta_2,
\end{equation*}
which implies that
\begin{equation}
 \label{nonint2}
\codim(Z_2,|\beta|)\geq \min\{\beta_1\cdot\beta_2\}.
\end{equation}

If there are $(-1)$-curves among $C_1,\cdots, C_{m}$, then we may assume that $C_1,\cdots, C_l$ are $(-1)$-curves ($l\in\{1,\cdots,m\}$) and $C_j$ ($j>l$) are nef. Then $|C_1|,\cdots,|C_l|$ are reduced points. The locus $Z_3$ of such $C$ is contained in the image of $|C_1|\times\cdots\times|C_l|\times|\gamma|\cong|\gamma|$ in $|\beta|$, where $\gamma$ is a nef, effective (possibly zero) divisor in $|\beta|$ such that $n_1C_1+\cdots+n_lC_l\in|\beta-\gamma|$. There are finitely many nef, effective divisors $\gamma$ such that $\beta-\gamma$ is a chain of $(-1)$-curves (possibly with multiplicities). By Lemma \ref{effdiv}, when $\beta$ is nef, 
\begin{equation*}
    \begin{aligned}
        \dim|\beta|-\dim|\gamma|&=\frac{(\beta-\gamma)\beta+(\beta-\gamma)\gamma-K_S\cdot(\beta-\gamma)}{2}\\
        &\geq\frac{(\beta-\gamma)(\beta-K_S)}{2}\geq\frac{1}{2}\min\{(\beta-K_S) L:\text{ $L$ is a $(-1)$-curve}\},
    \end{aligned}
\end{equation*}
which implies that
\begin{equation}
    \label{nonint3}
    \codim(Z_3,|\beta|)\geq\frac{1}{2}\min\{(\beta-K_S) L:\text{ $L$ is a $(-1)$-curve}\}.
\end{equation}

\begin{pro}
For every $i\in\bZ_{\geq0}$ and an arbitrary ample divisor $\beta_0$ on $S$,   $\beta=n\beta_0$ ($n\in\bZ_{>0}$) satisfies  the   condition $(\mathrm{A}_i)$ defined in Definition \ref{Ai} when $n$ is sufficiently large.
\end{pro}
\begin{proof}
By Lemma \ref{effdiv} and the boundedness of the ideal sheaves defining length $i+1$ subschemes of $S$, $\beta$ satisfies the condition (1) of $(\mathrm{A}_i)$ when $n$ is sufficiently large. By Bertini's theorem \cite[Chapter \uppercase\expandafter{\romannumeral2}, Theorem 8.18, p. 179]{Hartshorne}, the condition (2) of $(\mathrm{A}_i)$ is also satisfied by $\beta$ when $n$ is sufficiently large. It remains to verify that so is the condition (3) of $(\mathrm{A}_i)$.

We are going to bound from below the codimension of $|\beta|\setminus U$ in $|\beta|$ for $\beta=n\beta_0$. 
Assume there are nonzero, nef, effective divisors $\beta_1$ and $\beta_2$ with $n\beta_0=\beta_1+\beta_2.$ 
Then
\begin{equation}
\label{eq:beta12}
    n^2\beta^2_0=\beta_1^2+\beta_2^2+2\beta_1\cdot\beta_2.
\end{equation}
By (\ref{nonint2}), we need to bound $\beta_1\cdot\beta_2$ from below.
    If $\beta_1^2=0$ or $\beta_2^2=0$, we may assume $\beta_1^2=0$ and it follows that
\begin{equation}
\label{twonef1}
\beta_1\cdot\beta_2=(n\beta_0-\beta_1)\beta_1=n\beta_0\cdot\beta_1\geq n.        
\end{equation}

If $\beta_1^2>0$ and $\beta_2^2>0$, then by the Hodge index theorem \cite[Theorem 1.6.1, p. 88]{Lazars04} and (\ref{eq:beta12}),  $$\beta_1\cdot\beta_2\geq\sqrt{\beta_1^2\beta_2^2}\geq\sqrt{n^2\beta_0^2-2\beta_1\cdot\beta_2-1},$$
which implies that
\begin{equation}
\label{twonef2}
    \beta_1\cdot\beta_2\geq n\sqrt{\beta_0^2}-1.
\end{equation}
    
By (\ref{nonint1}), (\ref{nonint2}), (\ref{twonef1}), (\ref{twonef2}) and (\ref{nonint3}), if $S$ is minimal (thus $Z_3=\emptyset$), then
    \begin{equation*}
2\operatorname{codim}(|\beta|\setminus U,|\beta|)\geq\min\left\{\frac{n\beta_0(n\beta_0-K_S)}{2},2n\sqrt{\beta_0^2}-2,2n\right\},
    \end{equation*}
otherwise 
 \begin{equation*}
2\operatorname{codim}(|\beta|\setminus U,|\beta|)\geq\min\left\{\frac{n\beta_0(n\beta_0-K_S)}{2},2n\sqrt{\beta_0^2}-2,n+1\right\}.
    \end{equation*}
Hence when $n$ is sufficiently large, the condition (3) of $(\mathrm{A}_i)$ is satisfied by $\beta$. 
\end{proof}

By the argument above, we can calculate the codimensions of the non-integral locus in the following special cases.

\begin{exmp}
\label{exmp:codim}
When $S=\bP^2$ and $\beta=\cO_{\bP^2}(d)$ ($d\in\bZ_{>0}$), 
 \begin{equation*}
\operatorname{codim}(|\beta|\setminus U,|\beta|)=d-1.
    \end{equation*}
When $S=\bP^1 \times \bP^1$ and $\beta =\cO_{\bP^1}(a_1)\boxtimes\cO_{\bP^1}(a_2)$ ($a_1,a_2\in \bZ_{>0}$), 
 \begin{equation*}
\operatorname{codim}(|\beta|\setminus U,|\beta|)=\min\{a_1,a_2\}.
    \end{equation*}
When $S=S_1$ is the blowup of a point on $\bP^2$ and $\beta =\cO_{S_1}(ah-de_1)$ ($a\in\bZ_{>0}$, $a>d$, $h$ is the pullback of the hyperplane class on $\bP^2$ and $e_1$ is the exceptional divisor),
 \begin{equation*}
\operatorname{codim}(|\beta|\setminus U,|\beta|)=\min\{a-d,d+1\}.
    \end{equation*}
\end{exmp}

In the following, we provide facts about the fiber dimensions of two morphisms of our interest, which will be used to prove the stabilization of $b_i(M_\beta)$.
\begin{pro}[{\cite[Proposition 2.6]{MS23}}, {\cite[Corollary 1.3]{Yuan23}}]
\label{eqdim:h}
When $\beta\cdot H$ and $\chi$ are coprime, the Hilbert-Chow morphism $h_{\beta,\chi}$ has fibers of the same dimension $p_a(\beta)$.
\end{pro}

\begin{pro}
\label{eqdim:cC}
    The morphism $\tilde{\pi}^{[n]}:\cC^{[n]}\to|\beta|$ has fibers of the same dimension $n$.
\end{pro}
\begin{proof}
    We need to prove that for every curve $C\in|\beta|$, the Hilbert scheme $C^{[n]}$ has dimension $n$. Since $C$ is a curve on $S$, there exists a $0$-dimensional closed subscheme $W\subset C$ such that 
    $$C^\circ:= C\setminus W \cong \Spec(\bC[x,y]/f(x,y))$$ for some (possibly reducible) polynomial $f(x,y)\in\bC[x,y]$. 
    Suppose the underlying set of $W$ consists of $s$ distinct points $x_1,\cdots,x_s\in C$ ($s\in\bZ_{>0}$).
    By \cite[Theorem 1.1]{Luan23}, every irreducible component of $(C^\circ)^{[n]}$ has dimension $n$.
    Note that $(C^\circ)^{[n]}$ is an open subscheme of $C^{[n]}$. Let $Y$ be an irreducible component of $C^{[n]}$. If $Y\cap(C^\circ)^{[n]}$ is nonempty, then it is an irreducible component of $(C^\circ)^{[n]}$ and therefore 
    $$\dim Y=\dim Y\cap(C^\circ)^{[n]} =n.$$
    
    If $Y\cap(C^\circ)^{[n]}=\emptyset$, then let $\tau_k:S^{[k]}\to \mathrm{Sym}^k(S)$ $(k\in\bZ_{>0})$ be the Hilbert-Chow morphism associated to $S^{[k]}$. Denote by $H_{x,k}$ the fiber $\tau_{k}^{-1}(kx)$ ($x$ is a closed point of $S$) and set $H_{x,0}:=\Spec\bC$. It is well-known that $\dim H_{x,k}=k-1$ (\cite[Théorème \uppercase\expandafter{\romannumeral5}.3.2]{Briancon77}).
    Since $W^{[k]}$ is a closed subscheme of $S^{[k]}$ and $Y$ can be stratified according to the length $k_j$ of the subscheme supported at $x_j$ $(j=1,\cdots,\cdots,s)$,
    \begin{equation*}
        \begin{aligned}
            \dim Y&\leq\max_{(k_1,\cdots,k_s)\in\Lambda_s}\dim(H_{x_1,k_1}\times\cdots H_{x_s,k_s}\times (C^\circ)^{[n-k_1-\cdots-k_s]})\\
            &=\max_{(k_1,\cdots,k_s)\in\Lambda_s}\{[k_1-1]_++\cdots+[k_s-1]_++(n-k_1-\cdots-k_s)\}\\
            &<n,
        \end{aligned}
    \end{equation*}
     where $[l]_+:=\max\{l,0\}$  ($l\in\bZ$) and 
     $$\Lambda_s=\left\{(n_1,\cdots,n_s)\in\bZ_{\geq0}^s:0<\sum_{j=1}^sn_j\leq n\right\}.$$
     We conclude that $\dim C^{[n]}=n$.
\end{proof}

\subsection{Proof of the main result}
We use $H^*(-)$ ({\it resp}. $H_c^*(-)$) to denote the singular cohomology ({\it resp}. singular cohomology with compact support) with $\bQ$-coefficients.
\begin{thm}
\label{thm:betti}
    For any given $i\in\bZ_{\geq0}$, $\chi\in\bZ$, every  effective divisor $\beta$ satisfying $(\mathrm{A}_i)$ as defined in Definition \ref{Ai} and each integer $k\leq i$, we have
    $$\dim \IH^k(M_{\beta,\chi})=\lim_{m\to\infty}b_k(S^{[m]}).$$
\end{thm}
\begin{proof} 
By Theorem \ref{thm:chiindep}, $\dim\IH^k(M_{\beta,\chi})=b_k(M_{\beta,1})$, so it suffices to consider the case when $\chi=1$. We write $M_\beta:=M_{\beta,1}$ and $h_\beta:=h_{\beta,1}$.
Since the odd Betti numbers of $M_\beta$ are zero by Proposition \ref{oddvanish} and so are those of $S^{[m]}$ ($m\in\bZ_{\geq0}$) by Göttsche's formula \cite[Theorem 0.1]{Got90}, we may assume $i$ is even. Let $\beta$ satisfy $(\mathrm{A}_i)$. We divide the proof into the following steps, of which the first two are inspired by the argument in \cite[Section 5]{Sac19}.

{\it Step 1}. Compute $b_i(\cC^{[i+1]})$.
By the condition (1) of $(\mathrm{A}_i)$ and Proposition \ref{iveryamp}, $\cC^{[i+1]}$ is a projective bundle over $S^{[i+1]}$ with fibers of dimension
$$\dim|\beta|+(i+1)-2(i+1)\geq\frac{i}{2}.$$
Therefore, by the Leray-Hirsch Theorem \cite[Theorem 4D.1, p. 432]{Hatcher02} and Göttsche's formula,
\begin{equation}
\label{eq:pbdl}
b_i(\cC^{[i+1]})=\sum_{j\leq i/2}b_{2j}(S^{[i+1]})=\sum_{k\leq i}\lim_{m\to\infty}b_k(S^{[m]}).
\end{equation}

{\it Step 2.} Relate $b_i(\cC^{[i+1]})$ to $b_i(\cC_{U}^{[i+1]})$ with the notation as in (\ref{notation:cun}).
Denote by 
$$d(i,\beta):=\dim|\beta|+i+1$$ the dimension of $\cC^{[i+1]}$.
The pair $(\cC^{[i+1]},\cC^{[i+1]}\setminus \cC_{U}^{[i+1]})$ induces a long exact sequence (for example see \cite[Appendix B.2.6, p. 420]{PS08})
\begin{equation*}
    \cdots\to H^{l-1}(\cC^{[i+1]}\setminus \cC_{U}^{[i+1]}) \to H^l_c(\cC_{U}^{[i+1]})\to H^l(\cC^{[i+1]})\to H^l(\cC^{[i+1]}\setminus \cC_{U}^{[i+1]})\to \cdots.
\end{equation*}
When $l=2d(i,\beta)-i$, the condition (3) of $(\mathrm{A}_i)$ and Proposition \ref{eqdim:cC} imply that 
$$l>l-1>2\dim(\cC^{[i+1]}\setminus \cC_{U}^{[i+1]})$$
and therefore the long exact sequence above yields
\begin{equation}
\label{isom:ub}
H_c^{2d(i,\beta)-i}(\cC_{U}^{[i+1]})\cong H^{2d(i,\beta)-i}(\cC^{[i+1]}).
\end{equation}
Hence by Poincaré duality and (\ref{isom:ub}),
\begin{equation}
\label{biub}
    b_i(\cC^{[i+1]})=b_{2d(i,\beta)-i}(\cC^{[i+1]})=\dim H_c^{2d(i,\beta)-i}(\cC_{U}^{[i+1]})=b_i(\cC_{U}^{[i+1]}),
\end{equation}
where the last identity follows from \cite[Corollary B.25, p. 423]{PS08} and the smoothness of $\cC_{U}^{[i+1]}$.

{\it Step 2'.} Replace $(\cC^{[i+1]},\cC^{[i+1]}\setminus \cC_{U}^{[i+1]})$ by $(M_{\beta},M_{\beta}\setminus h_{\beta}^{-1}(U))$ in Step 1. By the condition (3) of $(\mathrm{A}_i)$, Proposition \ref{eqdim:h} and the argument of Step 1, it follows that for $k\leq i$,
\begin{equation}
\label{eq:bkmb}
  b_k(M_{\beta})=b_{2(\beta^2+1)-k}(M_{\beta})=\dim H_c^{2(\beta^2+1)-k}(h_{\beta}^{-1}(U))=b_k(h_{\beta}^{-1}(U)). 
\end{equation}

{\it Step 3.} Relate $b_i(\cC_{U}^{[i+1]})$ to $b_i(\cC_{U}^{[i+1+Nd_0]})$ for $N\in\bZ_{>0}$, where $d_0$ is the degree of the multisection $D_0$ of $\pi_0:\cC_{U_0}\to U_0$ in (\ref{multisec}). By Proposition \ref{n+d0},
\begin{equation}
\label{eq:indepn} b_i(\cC_{U}^{[i+1]})=b_i(\cC_{U}^{[i+1+d_0]})=\cdots=b_i(\cC_{U}^{[i+1+Nd_0]}).
\end{equation}

{\it Step 4.} Use the Abel-Jacobi map $\mathrm{AJ}:\cC_{U}^{[i+1+Nd_0]}\to \overline{J}^{i+1+Nd_0}_{\pi}= h_{\beta,\chi}^{-1}(U)$, where
$$\chi=(i+1+Nd_0)+1-p_a(\beta).$$ 
Since curves on $S$ are Gorenstein, such a map exists by \cite[(8.2.3)]{AK80}. It sends a length $i+1+Nd_0$ closed subscheme $Z\subset C$ ($C\in|\beta|$ integral) defined by the ideal sheaf $\cI_Z$ to the dual $\cI_Z^\vee$.
If $i+1+Nd_0> 2p_a(\beta)-2$, then 
$$h^1(F)=\dim\mathrm{Hom}(F,\omega_{C})=0,\quad h^0(F)=(i+1+Nd_0)+1-p_a(\beta)$$
by Serre duality and the Riemann-Roch formula, where $F\in h_{\beta,\chi}^{-1}(C)$ and $\omega_C$ is the dualizing sheaf of $C$.
Then
$\mathrm{AJ}^{-1}(F)\cong\bP(H^0(F)^\vee)$ ({\it cf}. \cite[Proposition B.5]{PT10}) and $\mathrm{AJ}$ is a $\bP^{i+1+Nd_0-p_a(\beta)}$-bundle map. Thus when $N$ is large enough,
\begin{equation}
  \label{eq:ajcor}
\begin{aligned}
b_i(\cC_{U}^{[i+1+Nd_0]})&=\sum_{j\leq i/2}b_{2j}(\overline{J}^{i+1+Nd_0}_{\pi}) &\text{(by Leray-Hirsch)}\\
&=\sum_{j\leq i/2} b_{2j}(\overline{J}^{p_a(\beta)}_{\pi})=\sum_{j\leq i/2}b_{2j} (h_{\beta}^{-1}(M_\beta)) &\\
    &=\sum_{k\leq i}b_k(M_{\beta}) &\text{(by (\ref{eq:bkmb}))},\\
\end{aligned}
\end{equation}
where the second equality follows by restricting (\ref{chiindep}) to $U$ and taking hypercohomologies.

{\it Step 5.} Conclude the result. 
By (\ref{eq:pbdl}), (\ref{biub}), (\ref{eq:indepn}) and (\ref{eq:ajcor}),
\begin{equation}
\label{eq:otoh}
    b_i(\cC^{[i+1]})=\sum_{k\leq i}b_k(M_\beta)=\sum_{k\leq i}\lim_{m\to\infty}b_k(S^{[m]}).
\end{equation}
If $i\geq2$, then $\beta$ also satisfies $(\mathrm{A}_{i-2})$. We can replace $i$ by $i-2$ in (\ref{eq:otoh}) and solve inductively that
$$b_k(M_\beta)=\lim_{m\to\infty}b_k(S^{[m]}),$$
which completes the proof.
\end{proof}

\section{Families of one-dimensional sheaves} 
\label{sec:pic}
In this section, $(S,H)$ is an arbitrary polarized surface.
Assume $\beta$ satisfies $(\mathrm{A}_0)$ defined in \ref{Ai}.
We will construct families of sheaves in $M_{\beta}$ over curves and use these families as testing curves to show the linear independence of some elements in 
$$\mathrm{Num}_{\bQ}(M_{\beta}):=(\Pic(M_{\beta})/\equiv_{\mathrm{num}})\otimes_\bZ\bQ,$$ 
where $\equiv_{\mathrm{num}}$ is the numerical equivalence. 
The Picard number $\rho(M_{\beta})$ of $M_{\beta}$ is defined to be the dimension of $\mathrm{Num}_{\bQ}(M_{\beta})$ as a $\bQ$-vector space ({\it cf}. \cite[Definition 1.1.7, p. 18]{Lazars04}).
There are three types of families that we are going to construct.

\subsection{Moving support}\label{moving curve}
Suppose there is a pencil 
$$\bP^1\subset|\beta|$$ such that $R\cap U_0\not=\emptyset$, where $U_0\subset |\beta|$ is the locus of smooth curves. There is a multisection $R$ of $h_1: h_{\beta}^{-1}(\bP^1)\rightarrow \bP^1$ whose degree is $\deg\tilde{\tau}\in\bZ_{>0}$. Let $\tilde{R}$ be the normalization of $R$ and let \[ \tilde{\tau}: \tilde{R} \rightarrow M_\beta \] be the composition of the normalization map and the closed embedding $R \hookrightarrow M_\beta$.
Recall that $\lambda([\cO_x])=h_{\beta}^*\cO_{|\beta|}(1)$ by Corollary \ref{cor:lambda0}. Then we have
\begin{equation}
    \label{degmov}   \deg  \tilde{\tau}^{\,*}\lambda([\cO_x])=\deg \tilde{\tau}>0.
\end{equation}

\subsection{Fixed smooth support}\label{smoothfiber}
In this section, we fix a smooth and connected curve $C_0\in\left|\beta\right|$. Let $\nu:C_0\times C_0\hookrightarrow C_0\times S$ be the inclusion map and let $y\in C_0$ be a closed point (also viewed as a closed point of $S$).
Let 
$$\cF_{C_0}:=\nu_*(\cO_{C_0\times C_0}((p_a(\beta)-1)\,C_0\times\{y\}+\D_0)),$$
where $\D_0$ is the diagonal of $C_0\times C_0$. Then $\cF_{C_0}$ is a $C_0$-flat family of sheaves in $h_\beta^{-1}(C_0)$. This is a direct generalization of the construction in \cite[\S 2.2, p. 491]{CC15}.

To compute the Chern character $\mathrm{ch}(\cF_{C_0})$ of $\cF_{C_0}$, we apply the Grothendieck-Riemann-Roch formula 
and obtain
\begin{equation*}
    \begin{aligned}
        \mathrm{ch}(\cF_{C_0})\mathrm{td}(C_0\times S)&=\nu_*[\mathrm{ch}(\cO_{C_0\times C_0}((p_a(\beta)-1)\,C_0\times\{y\}+\D_0))\mathrm{td}(C_0\times C_0)]\\
        &= \nu_*[(1+(p_a(\beta)-1)\,C_0\times\{y\}+\D_0)\mathrm{td}(C_0\times C_0)],
    \end{aligned}
\end{equation*}
where the second equality follows from $\D_0^2=(2-2p_a(\beta))\{y \}\times\{y \}$ (for simplicity of notation we do not distinguish different zero-cycles).
Therefore,
\begin{equation}
\label{grrc}
    \begin{aligned}
        \mathrm{ch}(\cF_{C_0})p_S^*\mathrm{td}(S)&= \nu_*[(1+(p_a(\beta)-1)\,C_0\times\{y\}+\D_0)p_2^*\mathrm{td}(C_0)]\\
        &= C_0\times C_0+\nu_*(\D_0)+(1-p_a(\beta))\{y\}\times\{y\},
    \end{aligned}
\end{equation}
where $p_S:C_0\times S\to S$ and $p_2:C_0\times C_0\to C_0$ are the projections to the second factors.

\subsection{Fixed singular support}\label{singularfiber}
Suppose $D=C_1+C_2\in\left|\beta\right|$, where $C_1$ and $C_2$ are two smooth, connected curves on $S$. Assume $C_2\cdot H\geq C_1\cdot H$ and that $C_1$ intersects $C_2$ transversally at $$N:=C_1\cdot C_2\geq 2$$ points, say at $y_1,\cdots,y_N$. Then there is a natural exact sequence of sheaves on $D$
\begin{equation*}
    0 \rightarrow \cO_D \rightarrow \cO_{C_1} \bigoplus  \cO_{C_2} \rightarrow \mathop{\bigoplus}_{j=1}^N \bC_{y_j} \rightarrow 0.
\end{equation*}
Denote by $g_i$ ($i=1,2$) the genus of $C_i$.
Then by the adjunction formula,
$$p_a(\beta)=g_1+g_2+N-1.$$
It is well-known  that there is a bijection between pure $1$-dimensional sheaves on $D$ and the triple $(\cF_1,\cF_2, \varphi_j: \cF_1|_{y_j} \rightarrow  \cF_2|_{y_j})$ where $\cF_i$ is a vector bundle on $C_i$ $(i=1,2)$ and $\varphi_j$ is a linear map (for example see \cite[Lemma 2.3]{NS97}). Now we use this bijection to construct a stable one-dimensional sheaf $\cL$ supported on $D$.

Let $\cL_1$ ({\it resp.} $\cL_2$) be a line bundle on $C_1$ ({\it resp.} $C_2$) of degree
$g_1+1$ ({\it resp.} $g_2+N-1$). 
By the Riemann-Roch formula,
\begin{equation*}
\label{chiL12}
    \chi(\cL_1)=\deg(\cL_1)+1-g_1=2\quad\textnormal{and}\quad\chi(\cL_2)=\deg(\cL_2)+1-g_2=N.
\end{equation*}
We can glue $\cL_i$ to obtain a line bundle $\cL$ on the nodal curve $D$ with 
\begin{equation*}\label{chiL}
\chi(\cL)=\chi(\cL_1)+\chi(\cL_2)-N=2
\end{equation*} 
as follows. At $y_j$ ($j=1,\cdots,N$),
we fix the isomorphisms $\cL_1|_{y_j}\cong\cL_2|_{y_j}\cong \bC_{y_j}$. Let $\cL$ be the kernel of the following surjection induced by evaluation maps
$$\pi:\iota_{1*}\cL_1\bigoplus \iota_{2*}\cL_2 \to \mathop{\bigoplus}_{j=1}^N \bC_{y_j},$$
where $\iota_i:C_i\to D$ is the inclusion map.

\begin{lem}
\label{lem:Lst}
    The torsion sheaf  $\iota'_{\ast} \cL$ on $S$ is stable with respect to $H$, where $\iota':D \hookrightarrow S$ is the closed embedding.
\end{lem}
\begin{proof}
    Since $\cL|_{C_i}\cong\cL_i$ by the construction of $\cL$, we have short exact sequences 
    $$0\lra\cK_2\lra\cL\stackrel{r_1}\lra\cL_1\lra0$$
    and
    $$0\lra\cK_1\lra\cL\stackrel{r_2}\lra\cL_2\lra0.$$
    By the construction of $\cL$, the slopes of the above sheaves are
    \begin{equation*}
    \begin{aligned}
        \mu(\cL_1)&=\frac{2}{C_1\cdot H},\ \ &\mu(\cK_2)&=\frac{\chi(\cL)-\chi(\cL_1)}{C_2\cdot H}=0,\\
        \mu(\cL_2)&=\frac{N}{C_2\cdot H},\ \ &\mu(\cK_1)&=\frac{\chi(\cL)-\chi(\cL_2)}{C_1\cdot H}=\frac{2-N}{C_1\cdot H}.
    \end{aligned}
    \end{equation*}
    Let $\cF$ be a nonzero proper subsheaf of $\cL$. If the support $\Supp(\cF)$ of $\cF$ is $D$, then 
    $$\mu(\cF)=\frac{\chi(\cF)}{\beta\cdot H}<\frac{\chi(\cL)}{\beta\cdot H}=\mu(\cL).$$
    Now assume $\mathrm{Supp}(\cF)=C_1$ or $C_2$.
    If the restriction
    $r_1|_{\cF}:\cF\to\cL_1$ is zero, then $\cF$ is a subsheaf of $\cK_2$. Since $\cK_2$ is a line bundle supported on $C_2$,
    $$\mu(\cF)\leq\mu(\cK_2)=0<\mu(\cL)=\frac{2}{\beta\cdot H}.$$
    If $r_1|_{\cF}$ is nonzero, then $\mathrm{Supp}(\cF)=C_1$ and $r_2|_{\cF}=0$. Hence $\cF$ is a subsheaf of $\cK_1$ and 
    $$\mu(\cF)\leq\mu(\cK_1)<\mu(\cL)$$
    since $N\geq 2$ by our assumption.
\end{proof}

Next, we construct a $C_1$-flat family of sheaves in $h_\beta^{-1}(D)$.
Let $\cE_0$ be the pull back of $\iota'_{*}\cL$ along the projection $p_S:C_1\times S\to S$. Then there is a short exact sequence of sheaves on $C_1\times S$
$$0\lra\cK\to\cE_0\lra \nu_{1*}(\cE_0|_{C_1\times C_1})\lra0,$$
where $\nu_1:C_1\times C_1\to C_1\times S$ is the inclusion.
Note that $\cE_0|_{C_1\times C_1}\cong p_2^*\cL_1$
is a line bundle, where $p_2:C_1\times C_1\to C_1$ is the projection to the second factor.

Write 
$$\cG_0:=(\cE_0|_{C_1\times C_1})\otimes \cO_{C_1\times C_1}(-\D_{1})\subset\cE_0|_{C_1\times C_1},$$
where $\D_{1}$ is the diagonal in $C_1\times C_1$.
Let $\cF_{C_1}$ be the subsheaf of $\cE_0\oplus \nu_{1*}\cG_0$ defined by
$$\cF_{C_1}(V)=\{(s,t)\in\cE_0(V)\oplus \cG_0(\nu_1^{-1}(V)):s|_{\nu_1^{-1}(V)}=t\}$$
for every open subset $V$ of $C_1\times S$.
Then
we have a commutative diagram of $C_1$-flat sheaves with exact rows
\begin{equation}
\label{diagram:ec1}
\begin{tikzcd}
{0}\arrow[r,""] & {\cK} \arrow[r,""] \arrow[d, "="]           & {\cF_{C_1}} \arrow[r,""] \arrow[d, hook, ""]  & {\nu_{1*}\cG_0}\arrow[d, hook, ""]\arrow[r,""] & {0}\\
{0}\arrow[r,""] &{\cK} \arrow[r, ""]  & {\cE_0} \arrow[r,""] & {\nu_{1*}(\cE_0|_{C_1\times C_1})} \arrow[r,""] & {0}.
\end{tikzcd}
\end{equation}

\begin{lem}
    The sheaf $\cF_{C_1}$ is a $C_1$-family of sheaves in $h_\beta^{-1}(D)$.
\end{lem}
\begin{proof}
    For every closed point $c\in C_1$, restricting the diagram (\ref{diagram:ec1}) of $C_1$-flat sheaves to $\{c\}\times S\cong S$ yields
    \begin{equation*}
\begin{tikzcd}
{0}\arrow[r,""] & {\cK_2} \arrow[r,""] \arrow[d, "="]           & {\cF_{C_1}|_{\{c\}\times S}} \arrow[r,"r_1'"] \arrow[d, hook, ""]  & {\cL_1(-c)}\arrow[d, hook, ""]\arrow[r,""] & {0}\\
{0}\arrow[r,""] &{\cK_2} \arrow[r, ""]  & {\cL} \arrow[r,"r_1"] & {\cL_1} \arrow[r,""] & {0}.
\end{tikzcd}
\end{equation*}
As a subsheaf of $\cL$, $\cF_{C_1}|_{\{c\}\times S}$ is pure and 
$$\chi(\cF_{C_1}|_{\{c\}\times S})=\chi(\cK_2)+\chi(\cL_1(-c))=\chi(\cL)-1=1.$$
It remains to show that $\cF_{C_1}|_{\{c\}\times S}$ is stable. Let $\cF\subset\cF_{C_1}|_{\{c\}\times S}$ be a nonzero proper subsheaf. As in the proof of Lemma \ref{lem:Lst}, we may assume $\mathrm{Supp}(\cF)$ is $C_1$ or $C_2$.
If the restriction $r_1'|_{\cF}:\cF\to\cL_1(-c)$ is zero, then $\cF$ is a subsheaf of $\cK_2$ and therefore
$$\mu(\cF)\leq\mu(\cK_2)=0<\mu(\cF_{C_1}|_{\{c\}\times S}).$$
If $r_1'|_{\cF}$ is nonzero, then $\mathrm{Supp}(\cF)=C_1$. We claim that in this case, $\mu(\cF)\leq0<\mu(\cF_{C_1}|_{\{c\}\times S})$ and thus $\cF_{C_1}|_{\{c\}\times S}$ is stable.
Otherwise  $\mu(\cF)>0$, then 
$$\mu(\cF)=\frac{\chi(\cF)}{C_1\cdot H}\geq\frac{1}{C_1\cdot H}\geq\frac{2}{\beta\cdot H}=\mu(\cL),$$
where the second inequality follows from our assumption that $C_2\cdot H\geq C_1\cdot H$. However, since $\cF$ is also a subsheaf of $\cL$, this
contradicts the stability of $\cL$ proved in Lemma \ref{lem:Lst}!
\end{proof}
\begin{rem}
    Indeed, our construction of the family $\cF_{C_1}$ can  be viewed as a family of Hecke modifications of $\cL$ at points on $C_1$. 
\end{rem}
Using the Grothendieck-Riemann-Roch formula, as in (\ref{grrc}) we have
\begin{equation}
\label{grrc1}
    \begin{aligned}
        &(\mathrm{ch}(\cF_{C_1})-\mathrm{ch}(\cE_0))p_S^*\mathrm{td}(S)\\
        =&(\mathrm{ch}(\nu_{1*}\cG_0)-\mathrm{ch}(\nu_{1*}(\cE_0|_{C_1\times C_1}))p_S^*\mathrm{td}(S)\\
        =&\nu_{1*}[p_2^*\mathrm{ch}(\cL_1)(\mathrm{ch}(\cO_{C_1\times C_1}(-\D_{1}))-1)p_2^*\mathrm{td}(C_1)]\\
        = &-\D_{1}-(g_1+1)\{y\}\times\{y\}.
    \end{aligned}
\end{equation}

\subsection{Intersection numbers}
\label{subsec:internum}
By the assertion 3 in \cite[Theorem 8.1.5, p. 216]{HL10}, if $C$ is a curve and $\cF_C$ is a $C$-flat family of sheaves in $M_\beta$ inducing a morphism $\phi_{\cF_C}:C\to M_\beta$, then
$$\phi_{\cF_C}^*\lambda(\mathfrak{a})=\det(p_{C!}(\cF_C\otimes p_S^*\mathfrak{a}))\quad (\mathfrak{a}\in K^0_\beta(S)),$$
where $p_C$ ({\it resp.} $p_S$) is the 
projection from $C\times S$ to $C$ ({\it resp.} $S$).
Thus the intersection number of $\lambda(\mathfrak{a})$ with $C$ is given by
\begin{equation}\label{deg}
\begin{aligned}
    &\deg(\phi_{\cF_C}^*\lambda(\mathfrak{a}))=\deg(\det(p_{C!}(\cF_C\otimes p_S^*\mathfrak{a})))\\
    =&\deg( c_1((p_{C!}(\cF_C\otimes p_S^*\mathfrak{a})))\\
    \stackrel{(\ast)}=&\deg (p_{C*}\{\mathrm{ch}(\cF_C)\cdot p_S^*\mathrm{ch}(\mathfrak{a})\cdot p_S^*\mathrm{td}(S)\}_3)\\
    =&\deg \left[ \{\mathrm{ch}(\cF_C)\cdot p_S^*\mathrm{td}(S)\}_1 \cdot ch_2( p_S^*\mathfrak{a})+\right.\\
    & \left. \{\mathrm{ch}(\cF_C)\cdot p_S^*\mathrm{td}(S)\}_2\cdot ch_1( p_S^*\mathfrak{a}) +r(\mathfrak{a})  \{\mathrm{ch}(\cF_C)\cdot p_S^*\mathrm{td}(S)\}_3\right] .
\end{aligned}
\end{equation}
where $\{\}_i$ means taking the degree $i$ part of a Chow class, $r(\mathfrak{a})$ is the rank of $\mathfrak{a}$ and $(\ast)$ follows from the Grothendieck-Riemann-Roch formula.

Now we will use our construction of the testing curves and  the calculations in \S \ref{moving curve}-\S \ref{singularfiber} to compute the intersection numbers of line bundles with these testing curves. 

Let $x$ be a closed point of $S$ and let $L$ be a divisor on $S$. Then a direct calculation shows $[\cO_x]\in K_\beta^0$ and $[-(\beta\cdot L)\cO_S+\cO_L]\in K_\beta^0$.
Note that
\begin{equation}
\label{ch:kc}
    \mathrm{ch}([\cO_x])=[x] \quad\textnormal{and}\quad\mathrm{ch}([-(\beta\cdot L)\cO_S+\cO_L])=-\beta\cdot L+L-\frac{1}{2}L^2.
\end{equation}
\begin{enumerate}
    \item If $C$ is the curve $\tilde{R}$ constructed in \S \ref{moving curve}, then  by (\ref{degmov}), 
$$\deg\tilde{\phi}^*\lambda([\cO_x])=\deg\tilde{\tau}>0.$$
    \item If $C$ is the curve $C_0$ constructed in \S \ref{smoothfiber}, then 
by (\ref{grrc}), (\ref{deg}) and (\ref{ch:kc}),
\begin{equation*}
\begin{split}
&\deg\phi_{\cF_{C_0}}^*\lambda([\cO_x])=0,\\ 
&\deg\phi_{\cF_{C_0}}^*\lambda([-(\beta\cdot L)\cO_S+\cO_L])=p_a(\beta)\beta\cdot L.
        \end{split}
    \end{equation*}

    \item Let $C$ be the curve $C_1$ constructed in \S \ref{singularfiber}.
    Since $\cE_0=p_S^*(\iota'_*\cL)$, for each $\mathfrak{a}\in K^0_\beta(S)$,
    \begin{equation}
    \label{cE0}
        \{\mathrm{ch}(\cE_0)\cdot p_S^*\mathrm{ch}(\mathfrak{a})\cdot p_S^*\mathrm{td}(S)\}_3=\{p_{S}^*\mathrm{ch}(\iota'_*\cL)\cdot p_S^*\mathrm{ch}(\mathfrak{a})\cdot p_S^*\mathrm{td}(S)\}_3=0.
    \end{equation}
    By (\ref{grrc1}), (\ref{deg}), (\ref{ch:kc}) and (\ref{cE0}), we have
    \begin{equation*}
        \begin{split}
           &\deg\phi_{\cF_{C_1}}^*\lambda([\cO_x])=0 , \\
           &\deg\phi_{\cF_{C_1}}^*\lambda([-(\beta\cdot L)\cO_S+\cO_L])=-C_1\cdot L+(g_1+1)\beta\cdot L. 
        \end{split}
    \end{equation*}
\end{enumerate}

\subsection{A lower bound of Picard numbers}

Fix a $\bQ$-basis $[L_1],\cdots,[L_{\rho}]$ for 
$\mathrm{Num}_{\bQ}(S)$, 
where $L_i\in\Pic(S)$ $(i=1,\cdots, \rho:=\rho(S))$ are smooth, connected effective divisors.
Write
\begin{equation}
\label{def:lambda_i}
 \lambda_0:=\lambda([\cO_x]),\ \lambda_i:=\lambda([-(\beta\cdot L_i)\cO_S+\cO_{L_i}]).
\end{equation}

We define the condition $(\mathrm{P})$ as follows.

\begin{defn}
\label{P}
    We say a divisor $\beta$ on $S$ satisfies $(\mathrm{P})$ if
    \begin{enumerate}
    \item $\beta$ satisfies $(\mathrm{A}_0)$ as defined in Definition \ref{Ai};
    \item $p_a(\beta)>0$;
        \item for $i=1,\cdots,\rho$, there exists a smooth and connected curve $L_i'\in|\beta-L_i|$ such that $L_i$ meets $L_i'$ transversally at $N_i\geq2$ points.
    \end{enumerate}
\end{defn}


\begin{thm}\label{picardnum}
When $\beta$ satisfies $(\mathrm{P})$ as defined in Definition \ref{P}, we have the following relation of Picard numbers
    $$\rho(M_\beta)\geq\rho(S)+1.$$
\end{thm}
 \begin{proof} 
We claim that $[\lambda_0],[\lambda_1],\cdots,[\lambda_\rho]$ which are defined in (\ref{def:lambda_i}) are linearly independent in $\mathrm{Num}_{\bQ}(M_\beta)$, thus the result follows. 
We prove the claim by contradiction. If to the contrary that $[\lambda_0],[\lambda_1],\cdots,[\lambda_\rho]$ were linearly dependent, there would exist $a_0,a_1,\cdots,a_\rho\in\bQ$, not all zero, such that $$a_0\lambda_0+a_1\lambda_1+\cdots+ a_\rho\lambda_\rho\equiv_{\mathrm{num}}0.$$ 
Then for every curve $C$ and every morphism $\phi_{\cF_C}:C\to M_\beta$, 
\begin{equation}
    \label{singmatr}
    a_0\deg \phi_{\cF_C}^*\lambda_0+a_1\deg \phi_{\cF_C}^*\lambda_1+\cdots,a_\rho\deg \phi_{\cF_C}^*\lambda_\rho=0.
\end{equation}

Write $\beta=\sum_{i=1}^\rho d_iL_i$
    for $d_i\in\bQ$. We may assume $d_1\not=0$.
Let $L_2',\cdots,L_r'$ be given as in Definition \ref{P} (3) and let $D_j=L_j+L_j'$ ($j=2,\cdots, \rho$).
By the calculations (1)-(3) at the end of \S \ref{subsec:internum}, we have the following table, where each entry is the intersection number (\ref{deg}) of a line bundle (a term in the first row) with a curve (a term in the first column), $n_i\in\bZ$ and the $*$'s are unimportant numbers. 
\begin{equation*}
\renewcommand\arraystretch{1.3}
    \begin{array}{c|ccccc}
\cdot & \lambda_0  & \lambda_1& \lambda_2 & \cdots &   \lambda_\rho \\ \hline 
\tilde{R} & \deg\tilde{\tau} & * & * &  & *  \\
C_0 & 0& p_a(\beta)\beta\cdot L_1 & p_a(\beta)\beta\cdot L_2 & \cdots & p_a(\beta) \beta\cdot L_\rho \\
D_2 & 0&  L_2\cdot L_1+n_2\beta\cdot L_1 & L_2\cdot L_2+n_2\beta\cdot L_2 &  &  L_2\cdot L_\rho+n_2\beta\cdot L_\rho  \\
\vdots & &  & &  &   \\
D_{\rho} & 0 &  L_\rho\cdot L_1+n_\rho\beta\cdot L_1 & L_\rho\cdot L_2+n_\rho\beta\cdot L_2 & &  L_\rho\cdot L_\rho+n_\rho\beta\cdot L_\rho  \\
\end{array}
\end{equation*}
The determinant of the intersection matrix can be calculated as
\[
\begin{aligned}
    &\det
\begin{pmatrix}
  \deg\tilde{\tau} & * & * &  & *\\
  0& p_a(\beta)\beta\cdot L_1 & p_a(\beta)\beta\cdot L_2 & \cdots & p_a(\beta) \beta\cdot L_\rho\\
  0&  L_2\cdot L_1+n_2\beta\cdot L_1 & L_2\cdot L_2+n_2\beta\cdot L_2 &  &  L_2\cdot L_\rho+n_2\beta\cdot L_\rho\\
  \vdots & \vdots & \vdots & \vdots & \vdots\\
  0 &  L_\rho\cdot L_1+n_\rho\beta\cdot L_1 & L_\rho\cdot L_2+n_\rho\beta\cdot L_2 & &  L_\rho\cdot L_\rho+n_\rho\beta\cdot L_\rho
\end{pmatrix}\\
=&\deg\tilde{\tau}\cdot p_a(\beta)\det
\begin{pmatrix}
   \beta\cdot L_1 & \beta\cdot L_2 & \cdots &  \beta\cdot L_\rho\\
    L_2\cdot L_1 & L_2\cdot L_2 &  &  L_2\cdot L_\rho\\
   \vdots & \vdots & \vdots & \vdots\\
   L_\rho\cdot L_1 & L_\rho\cdot L_2 & &  L_\rho\cdot L_\rho
\end{pmatrix}\\
=& d_1\deg\tilde{\tau}\cdot p_a(\beta)\det(L_k\cdot L_l)_{1\leq k,l\leq \rho}\not=0,
\end{aligned}
\]
which implies that the intersection matrix is nonsingular, a contradiction to (\ref{singmatr})!
 \end{proof}   

\section{Applications}
\label{sec:app}
In what follows, $S$ will be a del Pezzo surface with a divisor $\beta$ satisfying $(\mathrm{A_2})$ defined in Definition \ref{Ai}. We assume that $\beta\cdot H$ and $\chi$ are coprime.
Denote by $\CH^\ast(-)$ the Chow ring with $\bQ$-coefficients and by $H^\ast(-)$ the singular cohomology ring with $\bQ$-coefficients.

By Theorem \ref{markman}, we can identify $\CH^\ast(M_{\beta,\chi})$ with $H^\ast(M_{\beta,\chi})$ via the cycle class map. Then by Theorem \ref{thm:betti}, $\CH^1(M_{\beta,\chi})=\Pic(M_{\beta,\chi})\otimes_\bZ\bQ$ is generated by $\rho+1$ elements $(\rho=\rho(S))$, say by $\alpha_0,\alpha_1,\cdots,\alpha_\rho$.

\subsection{Normalized tautological classes}
\label{tautcl}
Let $\cE$ be a universal sheaf on $M_{\beta,\chi}\times S$. We consider the twisted Chern character of $\cE$ by a class $\alpha \in \CH^\ast(M_{\beta,\chi}\times S)$
\[ \mathrm{ch}^\alpha(\cE):=\mathrm{ch}(\cE)\cdot \exp(\alpha)=\mathop{\sum} \limits_{i\ge 0} \mathrm{ch}^\alpha_i(\cE),  \]
where $\exp(\alpha):=\sum_{n\geq0} \frac{\alpha^n}{ n!}$ and $\ch_i^{\alpha}(\cE)\in\CH^i(M_{\beta,\chi}\times S)$ is the degree $i$ part of $\ch^\alpha(\cE)$. Recall that we have two natural projections $p: M_{\beta,\chi} \times S \rightarrow M_{\beta,\chi}$ and $q: M_{\beta,\chi}\times S \rightarrow S$.  Given a  class  $\gamma \in \CH^j(S)$, it will produce\begin{equation} \label{chernclass}
    p_\ast(\mathrm{ch}_k^\alpha(\cE) \cdot q^\ast \gamma) \in \CH^{j+k-2}(M_{\beta,\chi}).
\end{equation}  
As $\gamma$ runs over a basis for $\CH^\ast(S)$, the number of divisor classes obtained by (\ref{chernclass}) is \[1+b_2(S)+1=b_2(M_{\beta,\chi})+1\] by Theorem \ref{thm:betti}. This number combined with Theorem \ref{beauville} suggests that there is only one relation between these $b_2(M_{\beta,\chi})+1$ divisor classes. To obtain free generators for divisor classes, we follow the normalization procedure in \cite[\S 2.1]{PS23}. Since $H^1(S,\cO_S)=0$, we have ({\it cf}.~\cite[Chapter III, Excercise 12.6, p. 292]{Hartshorne})
    \begin{equation}\label{productpicard}
        \CH^1(M_{\beta,\chi}\times S) \cong  p^\ast \CH^1(M_{\beta,\chi})\oplus q^\ast \CH^1(S).
    \end{equation}

\begin{lem}\label{twist}
   Let $\cE$ be a universal sheaf, then there is a unique divisor class $\alpha_\cE \in \CH^1(M_{\beta,\chi} \times S)$ such that 
   \begin{enumerate}
       \item the following classes vanish
       \begin{equation}\label{unique}
       p_\ast (\ch_2^{\alpha_\cE}(\cE))=0\in\CH^0(M_{\beta,\chi}\times S),\quad    p_\ast (\ch_2^{\alpha_\cE}(\cE)\cdot q^\ast K_S)=0\in\CH^1(M_{\beta,\chi}\times S);
   \end{equation}
   \item the second factor of $\alpha_\cE$ under the identification (\ref{productpicard}) lies in $\bQ [q^\ast K_S]\subset q^\ast\CH^1(S)$.
   \end{enumerate}
\end{lem}
\begin{proof}
By the condition (2) above,  we can write $$\alpha_\cE=u_0p^*\alpha_0+u_1p^*\alpha_1+\cdots+ u_\rho p^*\alpha_\rho+v q^* K_S$$ 
for $u_l,v\in \bQ$ $(l=0,\cdots,\rho)$ to be determined.
By Proposition \ref{firstchern} and Corollary \ref{cor:lambda0}, 
\[\ch_2^{\alpha_\cE}(\cE)=\ch_2(\cE)+ \alpha_\cE\cdot c_1(\cE)=\ch_2(\cE)+\alpha_\cE\cdot(p^*\lambda_0+q^*\beta)\in \CH^2(M_{\beta,\chi} \times S).
 \]
Hence the two equations in (\ref{unique}) correspond respectively to
\begin{equation*}
\begin{cases}
    (K_S\cdot\beta)v=-p_*(\ch_2(\cE))\\
    (K_S^2v+(K_S\cdot\beta)u_0)\alpha_0+(K_S\cdot\beta)(u_1\alpha_1+\cdots+ u_\rho \alpha_\rho)=-p_*(\ch_2(\cE)\cdot q^*K_S).
\end{cases}
\end{equation*}
Since $K_S^2>0$ and $K_S\cdot\beta<0$, $u_l$ and $v$ are uniquely determined.
\end{proof}

Now we are able to define the tautological classes for $M_{\beta,\chi}$ to generalize the construction in \cite{PS23}.
We introduce the notation for generators of $\Pic(S)$ as follows. When $S=\bP^2$, $\Pic(S)$ is freely generated by $h:=\cO_{\bP^2}(1)$. When $S=\bP^1\times\bP^1$, $\Pic(S)$ is freely generated by $h_1:=\cO_{\bP^1}(1)\boxtimes\cO_{\bP^1}$ and $h_2:=\cO_{\bP^1}\boxtimes\cO_{\bP^1}(1)$. When $S=S_\delta$ is the blowup of $\delta$ points on $\bP^2$ ($\delta=1,\cdots,8$), we denote by $h$ the pullback of $\cO_{\bP^2}(1)$ to $S_\delta$ and by $e_1,\cdots,e_\delta$ the exceptional divisors. The Picard group $\Pic(S_\delta)$ is freely generated by $h,e_1,\cdots,e_{\delta}$. 
Let $L_1,\cdots,L_\rho$ be the generators of $\Pic(S)$ as above.
\begin{defn}
\label{def:taut}
   Let $\alpha_\cE$ be the unique class in Lemma \ref{twist}. For $k\in\bZ_{\geq0}$, define the (normalized) tautological classes to be the following
   \[c_k(0):=p_\ast( \ch_{k+1}^{\alpha_\cE}(\cE))\ \in\ \CH^{k-1}(M_{\beta,\chi}),\quad c_k(2):=p_\ast( \ch_{k+1}^{\alpha_\cE}(\cE)\cdot q^\ast[\mathrm{pt}])\ \in\ \CH^{k+1}(M_{\beta,\chi}), \]
   \[c_k(1,j):=p_\ast( \ch_{k+1}^{\alpha_\cE}(\cE)\cdot q^*L_j)\ \in \CH^{k}(M_{\beta,\chi})\quad (j=1,\cdots,\rho),\]
 where $\mathrm{pt}$ denotes a closed point of $S$.
\end{defn}   

Note that by (\ref{unique}), $c_1(0)=0$ and there is a relation among $\{c_1(1,j)\}_{j=1}^\rho$. The following lemma shows that the tautological classes are well-defined.

\begin{lem}
   The tautological classes in Lemma \ref{twist} are independent of the choice of a universal sheaf $\cE$. 
\end{lem}
\begin{proof}
Let $\cE$ and $\cE'$ be two universal sheaves on $M_{\beta,\chi} \times S$. Then there is a line bundle $B$ on $M_{\beta,\chi}$ such that $\cE'=\cE \otimes p^\ast B$. Let $\alpha':=\alpha_\cE-p^*B$. Then
$$\ch^{\alpha'}(\cE')=\ch(\cE')\cdot\exp(\alpha')=(\ch(\cE)\cdot p^*\ch(B))\cdot(\exp(\alpha_{\cE})\cdot p^*\ch(B)^{-1})=\ch^{\alpha_{\cE}}(\cE).$$
In particular, $\ch_2^{\alpha'}(\cE')=\ch_2^{\alpha_{\cE}}(\cE)$. It follows from the uniqueness of $\alpha_{\cE'}$ in Lemma \ref{twist} that $\alpha'=\alpha_{\cE'}$.
Thus for every $k\in\bZ_{\geq0}$,
$$\ch_{k+1}^{\alpha_\cE}(\cE)=\ch_{k+1}^{\alpha_{\cE'}}(\cE')$$
and the result follows.
\end{proof}

\begin{pro}
    The tautological classes $c_k(0),c_k(2),c_k(1,j)$ ($k\in\bZ_{\geq0},j=1,\cdots,\rho$) generate the cohomology ring $H^\ast(M_{\beta,\chi})$ as a $\bQ$-algebra. 
\end{pro}
\begin{proof}
    Since $1\in H^0(S), L_j\in H^2(S)$ and the Poincaré dual of $\mathrm{pt}$ form a basis for $H^\ast(S)$, the result follows from Theorem \ref{beauville}.
    \end{proof}

\begin{cor}
\label{cor:biineq}
    For any $i\in\bZ_{\geq0}$, we have
    $$b_{2i+4}(M_{\beta,\chi})\leq \lim_{n\to\infty}b_{2i+4}(S^{[n]}).$$
\end{cor}
\begin{proof}
 The proof resembles that of \cite[Theorem 1.2 (b)]{PS23}. By Göttsche's formula \cite[Theorem 0.1]{Got90}, the Poincaré polynomial $p(S^{[n]},z):=\sum_{l\geq0}b_l(S^{[n]})z^l$ of $S^{[n]}$ ($n\in\bZ_{\geq0}$) can be expressed as 
    \begin{equation*}
    \begin{aligned}
         \sum_{n=0}^\infty p(S^{[n]},z)t^n&=\prod_{m=1}^\infty\frac{(1+z^{2m-1}t^m)^{b_1(S)}(1+z^{2m+1}t^m)^{b_1(S)}}{(1-z^{2m-2}t^m)^{b_0(S)}(1-z^{2m}t^m)^{b_2(S)}(1-z^{2m+2}t^m)^{b_0(S)}}\\
         &=\prod_{m=1}^\infty(1-z^{2m-2}t^m)^{-1}(1-z^{2m}t^m)^{-b_2(S)}(1-z^{2m+2}t^m)^{-1}\\
        &=\prod_{m=1}^\infty\left(\sum_{l=0}^\infty z^{(2m-2)l}t^{ml}\right)\left(\sum_{l=0}^\infty z^{2ml}t^{ml}\right)^{b_2(S)}\left(\sum_{l=0}^\infty z^{(2m+2)l}t^{ml}\right)\\
        &=(1+t+t^2+\cdots)(1+z^2t+z^4t^2+\cdots)^{b_2(S)}(1+z^4t+z^8t^2+\cdots)\\
        &\times\prod_{m=2}^\infty\left(\sum_{l=0}^\infty z^{(2m-2)l}t^{ml}\right)\left(\sum_{l=0}^\infty z^{2ml}t^{ml}\right)^{b_2(S)}\left(\sum_{l=0}^\infty z^{(2m+2)l}t^{ml}\right).
    \end{aligned}
    \end{equation*}
By this expression, for $m\geq2i+4$,
$$b_{2i+4}(S^{[m]})=\lim_{n\to\infty}b_{2i+4}(S^{[n]}),$$
which equals the coefficient of $z^{2i+4}$ in 
\begin{equation}
\label{eq:zseries}
      \begin{aligned}
        &(1+z^2+z^4+\cdots)^{b_2(S)}(1+z^4+z^8+\cdots)\\
        &\times\prod_{m=2}^\infty\left(\sum_{l=0}^\infty z^{(2m-2)l}\right)\left(\sum_{l=0}^\infty z^{2ml}\right)^{b_2(S)}\left(\sum_{l=0}^\infty z^{(2m+2)l}\right).
    \end{aligned}
\end{equation}
On the other hand, the number of generators in $H^{2i+4}(M_\beta)$ given by tautological generators equals the number of monomials of degree $2i+4$ in the following generating series
\begin{equation}
\label{eq:tautseries}
      \begin{aligned}
        &\left[\left(\sum_{l=0}^\infty c_0(2)^l\right)\prod_{j=2}^{b_2(S)}\left(\sum_{l=0}^\infty c_1(1,j)^l\right)\right]\left(\sum_{l=0}^\infty c_1(2)^l\right)\\
        &\times\prod_{m=2}^\infty\left(\sum_{l=0}^\infty c_m(0)^l\right)\left[\prod_{j=1}^{b_2(S)}\left(\sum_{l=0}^\infty c_m(1,j)^l\right)\right]\left(\sum_{l=0}^\infty c_m(2)^l\right).
    \end{aligned}
\end{equation}
By comparing (\ref{eq:zseries}) and (\ref{eq:tautseries}), it is immediate to see that there are $\lim_{n\to\infty}b_{2i+2}(S^{[n]})$ many generators in $H^{2i+2}(M_{\beta,\chi})$, which are possibly linearly dependent. Thus the result follows.
\end{proof}

\begin{rem}
\label{betterai}
The result of Theorem \ref{thm:betti} is still true if the following conditions hold
\begin{enumerate}
\item $\beta$ is $\min\{i,2\}$-very ample and
$$2\dim|\beta|\geq 3i+2;$$
\item a general curve in $|\beta|$ is smooth and connected;
\item the codimension $\operatorname{codim}(|\beta|\setminus U,|\beta|)$ of the locus of non-integral curves in $|\beta|$ satisfies 
$$2\operatorname{codim}(|\beta|\setminus U,|\beta|)>i+1.$$
\end{enumerate}
Indeed, by the proof of Proposition \ref{iveryamp}, even if $\beta$ is not $i$-very ample, there is an open subset $V\subset S^{[i+1]}$, such that $\sigma_i^\circ:=\sigma_i|_{\sigma_i^{-1}(V)}:\sigma_i^{-1}(V)\to V$ is a projective bundle with fibers of dimension
$f\geq i/2.$ Since $\sigma_i^{-1}(V)\cap\cC_U^{[i+1]}\not=\emptyset$, we can take the unique irreducible component $Y$ of $\cC^{[i+1]}$ that intersects $\sigma_i^{-1}(V)\cup\cC_U^{[i+1]}$. Applying the decomposition theorem \cite[Théorème 6.2.5, p. 163]{BBD82} to $\sigma_i|_{Y}:Y\to S^{[i+1]}$, we have
    $$\begin{aligned}
    R(\sigma_i|_{Y})_*\mathrm{IC}_{Y}&=\cP\oplus\bigoplus_{l=-f}^{f}\mathrm{IC}(S^{[i+1]},R^{l+f}\sigma^\circ_{i*}\bQ)[-l]\\
    &=\cP\oplus\bigoplus_{m=0}^f\bQ[f-2m+2i+2],
    \end{aligned}$$
    where $\cP$ is a summand of $R(\sigma_i|_{Y})_*\mathrm{IC}_{Y}$.
Applying $\bH^{i-f-(2i+2)}(-)$ to $R(\sigma_i|_{Y})_*\mathrm{IC}_{Y}$, we have
    $$\dim\IH^i(Y)\geq\sum_{k\leq i}\lim_{n\to\infty}b_k(S^{[n]}).$$
    By \cite[Theorem 6.7.4, p. 113]{Maxim19} (the intersection cohomology version of Step 2 in the proof of Theorem \ref{thm:betti}),
    $$\dim\IH^i(Y)= b_i(\cC_U^{[i+1]}).$$
Hence by the proof of Theorem \ref{thm:betti}, 
$$\sum_{k\leq i}b_k(M_\beta)\geq \sum_{k\leq i}\lim_{n\to\infty}b_k(S^{[n]})$$
and Corollary \ref{cor:biineq} forces $b_k(M_\beta)=\mathop{\lim} \limits_{n\to\infty}b_k(S^{[n]})$ for every integer $k\leq i$.
\end{rem}

\begin{exmp}
\label{exmp:special}
Let the notation be as in Example \ref{exmp:codim}.
    When $S=\bP^2$ and $\beta=\cO_{\bP^2}(d)$, it follows from Remark \ref{betterai}, Example \ref{exmp:codim} and \cite[Lemma 2.6]{CGKT20} that for every integer $i\leq 2d-4$,
    $$\dim\IH^i(M_{\beta,\chi})=\lim_{n\to\infty}b_i(S^{[n]}).$$
The same result holds when 
$$S=\bP^1\times\bP^1,\quad\beta=\cO_{\bP^1}(a_1)\boxtimes\cO_{\bP^1}(a_2),\quad i\leq2\min\{a_1,a_2\}-2$$ or
when 
$$S=S_1,\quad\beta=\cO_{S_1}(ah-de_1),\quad i\leq 2\min\{a-d,d+1\}-2.$$
\end{exmp}

\subsection{P=C conjecture on del Pezzo surfaces and enumerative geometry}
Now we use the normalized tautological classes in \S \ref{tautcl}  to formulate a $P=C$ conjecture for moduli spaces of one dimensional sheaves on del Pezzo surfaces.

\subsubsection{Chern filtration and Perverse filtration}
We use $c_k(\Vec{a})$ to denote the tautological classes defined in Definition \ref{def:taut}, where $\Vec{a}$ stands for $(0)$, $(2)$ or $(1,j)$ ($j=1,\cdots,\rho$). They will produce an increasing filtration 
$C_{\bullet} H^\ast(M_{\beta,\chi})$ given by
\begin{equation}
    C_i H^\ast(M_{\beta,\chi}):=\  \hbox{the subspace generated by all }\ \mathop{\prod}_{l=1}^m c_{k_l}(\Vec{a}_{l}) \ \hbox{with}\ \mathop{\sum} \limits_{l=1}^m k_l\le i
\end{equation}
as in \cite{KPS} and \cite{KLMP24} for $S=\bP^2$. The filtration $C_{\bullet} H^\ast(M_{\beta,\chi})$ is called the Chern filtration on $H^\ast(M_{\beta,\chi})$.

An important geometric feature of the moduli space $M_{\beta,\chi}$ is the proper Hilbert-Chow morphism $h_{\beta,\chi}: M_{\beta,\chi} \rightarrow |\beta|$, which is a weak abelian fibration. The perverse filtration 
$P_\bullet H^{\ast}(M_{\beta,\chi})$ is defined as
\begin{equation}
\label{def:pervfil}
    P_iH^m(M_{\beta,\chi}):=\im \big\{H^{m-\dim|\beta|}(|\beta|,\,^{\textbf{p}}\tau_{\leq i}(R(h_{\beta,\chi})_{*}\bQ_{M_{\beta,\chi}}[\dim|\beta|]))\to H^m(M_{\beta,\chi}) \big\}
\end{equation}
using the language of perverse sheaves, where $^{\textbf{p}}\tau_{\leq i}$ is the perverse truncation functor \cite[\S 2.1, p. 57]{BBD82}. Due to de Cataldo-Migliorini
\cite{dCM05}, there is a more explicit characterization of the perverse filtration on $H^\ast(M_{\beta,\chi})$ as follows. Let 
\[\eta:=\cup c_1(h_\beta^\ast \cO_{|\beta|}(1)):\ \  H^\ast(M_{\beta,\chi}) \rightarrow H^{\ast+2}(M_{\beta,\chi}) \]
be the ring homomorphism defined by the cup product with respect to $c_1(h_\beta^\ast \cO_{|\beta|}(1))$. Then \cite[Proposition 5.2.4]{dCM05} shows that the perverse filtration associated to $h_{\beta,\chi}: M_{\beta,\chi}\rightarrow |\beta|$ can be given as 
\begin{equation}\label{Pfiltra}
    P_i H^m(M_{\beta,\chi})=\mathop{\sum} \limits_{l\ge 1} \ker(\eta^{\dim |\beta|+l+i-m})\cap \im(\eta^{l-1}) \cap H^m(M_{\beta,\chi}).
 \end{equation}
 
\begin{conj}\label{P=C}
  With the notation as above, we have the following $P=C$ identity 
  \begin{equation*}
      P_\bullet H^{\ast}(M_{\beta,\chi})=C_\bullet H^{\ast}(M_{\beta,\chi}).
  \end{equation*}
\end{conj}
We summarize the previously known evidence for Conjecture \ref{P=C}.
\begin{enumerate}
    \item Maulik-Shen-Yin in \cite[Theorem 0.6]{MSY23}  generalized the Beauville decomposition for abelian schemes over a base to dualizable abelian fibrations and used this tool to prove $ C_\bullet H^{\le 2r-4}(M_{\beta,\chi}) \subset P_\bullet H^{\le 2r-4}(M_{\beta,\chi})$ for $S=\bP^2$ and $\beta=\cO_{\bP^2}(r)$ ($r\geq3$).
    \item Kononov-Pi-Shen \cite{KPS} verified the conjecture for  $S=\bP^2$ and $\beta=\cO_{\bP^2}(d)$ for $d\le 4$. Later Yuan \cite{Yuan23a} verified $ P_\bullet H^{\le 4}(M_{\beta,\chi})=C_\bullet H^{\le 4}(M_{\beta,\chi})$ for $S=\bP^2$ and $\beta=\cO_{\bP^2}(l)$ ($l\geq4$). 
\end{enumerate}
Now we give the first evidence on one side of Conjecture \ref{P=C} for any del Pezzo surface. This relies on the framework of Maulik-Shen-Yin \cite[Theorem 2.5]{MSY23} on the generalized Beauville decomposition for the dualizable abelian fibration.
Recall that $U \subset |\beta|$ is the open subset parametrizing integral curves in $|\beta|$ and $h_{\beta,\chi}^{-1}(U)$ is isomorphic to the relative compactified Jacobian $J:=\overline{J}^{\chi+p_a(\beta)-1}_{\pi}$.
\begin{pro}\label{CinP}
 The inclusion $C_\bullet H^{i}(M_{\beta,\chi}) \subset P_\bullet H^{i}(M_{\beta,\chi})$ holds for
 $$i\leq 2\codim(|\beta|\setminus U,|\beta|)-2.$$
\end{pro}
\begin{proof}
Applying the long exact sequence 
associated to the pair $(M_{\beta,\chi},h_{\beta,\chi}^{-1}(U))$ as in the proof of Theorem \ref{thm:betti}, we have 
 \begin{equation*}
     \cdots \rightarrow H^{l-1}(M_{\beta,\chi}\setminus h_{\beta,\chi}^{-1}(U)) \rightarrow H^l_c(h_{\beta,\chi}^{-1}(U))\rightarrow H^l(M_{\beta,\chi})  \rightarrow H^{l}(M_{\beta,\chi}\setminus h_{\beta,\chi}^{-1}(U)) \rightarrow \cdots
 \end{equation*}
 Combining with the dimension reason, we have \[H^{2\beta^2+2-l}(M_{\beta,\chi})\cong  H^l(M_{\beta,\chi})^\vee \cong  H^l_c(h_{\beta,\chi}^{-1}(U))^\vee\cong H^{2\beta^2+2-l}(h_{\beta,\chi}^{-1}(U))\] for $l >  2\dim (M_{\beta,\chi}\setminus h_{\beta,\chi}^{-1}(U))$. Thus the restriction map $H^i(M_{\beta,\chi}) \rightarrow H^i(h_{\beta,\chi}^{-1}(U))$ is an isomorphism for $i < 2 \codim (|\beta|\setminus U,|\beta|)-1$. 
The relative compactified Jacobian $J \rightarrow U$ is a dualizable abelian fibration in the sense of \cite[\S 1.4]{MSY23} and the argument in \cite[5.3.3]{MSY23}  carries over to show $C_\bullet H^{\ast}(h_{\beta,\chi}^{-1}(U)) \subset P_\bullet H^{\ast}(h_{\beta,\chi}^{-1}(U))$, from which the inclusion $C_\bullet H^{*}(M_{\beta,\chi}) \subset P_\bullet H^{*}(M_{\beta,\chi})$ follows.
\end{proof}

\begin{pro}\label{C=PH2}
The identity  $C_\bullet H^{\le 2}(M_{\beta,\chi}) =P_\bullet H^{\le 2}(M_{\beta,\chi})$ holds.  
\end{pro}
\begin{proof}
Let $N:=\dim |\beta|$. 
It follows from (\ref{chiindep}) that
$$R(h_{\beta,\chi})_*\bQ_{M_{\beta,\chi}}[N]\cong\bigoplus_{l=0}^{2p_a(\beta)}\mathrm{IC}\left(\wedge^lR^1\pi_{0*}\bQ_{\cC_{U_0}}\right)[-l].$$
Hence by the definition of the perverse filtration (\ref{def:pervfil}),
$$P_0H^2(M_{\beta,\chi})\cong H^{2-N}(|\beta|,\bQ[N])=H^2(|\beta|)$$
is one-dimensional. Indeed  by (\ref{Pfiltra}) and  Corollary \ref{cor:lambda0}, $P_0H^2(M_{\beta,\chi})=\bQ [c_0(2)]$.  By (\ref{Pfiltra}) again, we have $P_2H^2(M_{\beta,\chi})=H^2(M_{\beta,\chi})$ and 
\begin{equation*}
    \begin{split}
        P_1H^2(M_{\beta,\chi})=\ker(\eta^{N}: H^2(M_{\beta,\chi})\rightarrow H^{2N+2}(M_\beta))+\im(\eta).
    \end{split}
\end{equation*}
 Since $\eta^{N}$ is the fiber class of $h_{\beta,\chi}$ and $c_2(0)$ is relatively ample, $\eta^{N}.c_2(0) \neq 0$ and thus $c_2(0) \notin \ker(\eta^{N})$, which proves $c_2(0) \in P_2H^2(M_{\beta,\chi})\setminus P_1H^2(M_{\beta,\chi})$.  
Then the result follows from Proposition \ref{CinP}.    
\end{proof}

\subsubsection{Relation to enumerative geometry of local surfaces}
The total space $X=\mathrm{Tot}(K_S)$ of the canonical divisor $K_S$ is a non-compact Calabi-Yau $3$-fold for a del Pezzo surface $S$. The moduli space $M_{\beta,\chi}$ provides a mathematical definition of the moduli space of 2D-branes
on $X$, which is predicted from physics, see \cite{GV} \cite{HST}, \cite{KL12} and \cite{MT18}. There are two $\mathfrak{sl}_2$-actions on $M_{\beta,\chi}$ given by Lefschetz triples associated to two Kähler classes $c_2(0)$ and $c_0(2)$. Thus $ \IH^{\ast} (M_{\beta,\chi})$ is an $\mathfrak{sl}_2 \times  \mathfrak{sl}_2$ representation space and Conjecture  \ref{P=C} identifies the  irreducible representations and the graded pieces of the perverse filtration.  Then the refined BPS numbers can be mathematically defined as 
\begin{equation*}
    n_\beta^{i,j}:=\dim \gr^P_i \IH^{i+j}(M_{\beta,\chi})
\end{equation*}
which is independent of the choice of $\chi$ by Theorem \ref{thm:chiindep}.  A direct consequence in the proof of 
 Proposition \ref{C=PH2} is the following partial calculation  of refined BPS numbers.
\begin{cor}
With the notation as in  Proposition \ref{C=PH2}, we have
$$n_\beta^{0,2}=1,\ n_\beta^{1,1}=\rho(S)-1,\ n_\beta^{2,0}=1.$$
\end{cor}
We make the following stronger stabilization conjecture. 
\begin{conj}\label{conj:BPS}
The refined BPS number $ n_\beta^{i,j}$ stabilizes when $\beta$ is sufficiently positive.
\end{conj}
The stronger Conjecture \ref{conj:BPS}   implies  Conjecture \ref{conj} for del Pezzo surfaces.
For  $S=\bP^2$ and $\beta=\cO_{\bP^2}(d)$,  a conjectural product formula 
\begin{equation}\label{equ:bps}
    \mathop{\sum} \limits_{i,j} n_{\beta}^{i,j}q^it^j=\mathop{\prod} \limits_{i\ge 0 } \frac{1}{(1-(qt)^iq^2)(1-(qt)^iq^2t^2)(1-(qt)^it^2)},\ \ \ i+j\le 2d-4
\end{equation}
in \cite[Conjecture 0.1]{KPS} implies Conjecture \ref{conj:BPS}. It is natural to ask the following
\begin{question}
    Is there a similar formula for $\mathop{\sum} \limits_{i,j} n_{\beta}^{i,j}q^it^j$ as in (\ref{equ:bps}) for a del Pezzo surface?
\end{question}

\section{Reducible moduli spaces on certain surfaces of general type}
\label{sec:other}


In this section, we provide an example where the asymptotic irreducibility fails for certain surfaces of general type, which indicates more complicated behaviors may happen for surfaces of general type.
Let $S$ be a surface of general type, i.e., $K_S$ is big. 
We make the following assumption for $S$ and $\beta$.
\begin{assu}
\label{assu:irrd}
The following conditions hold:
    \begin{enumerate}
    \item a general curve in $|K_S| $ is smooth and connected;
    \item $q:=h^1(\cO_S) \ge 1$;
    \item $\beta_n\in \pic(S)$ such that $\beta_n \cong\cO_S(nK_S)$ for an integer $n\ge 2$.
\end{enumerate}
\end{assu}

\begin{pro}\label{counter} Under Assumption \ref{assu:irrd},
    the moduli space $M_{\beta_n}$ is not irreducible for all $n\ge 2$.
\end{pro}
\begin{proof}
By the Riemann-Roch formula and Serre duality,
\begin{equation}
\label{dimks}
    \begin{aligned}
        \dim|K_S|&=\chi(K_S)+h^1(K_S)-h^2(K_S)-1\\
        &=\left(\frac{K_S(K_S-K_S)}{2}+\chi(\cO_S)\right)+h^1(\cO_S)-h^0(\cO_S)-1\\
        &=\chi(\cO_S)+q-2.
    \end{aligned}
\end{equation}
For a general $C\in|K_S|$, which is smooth and connected by the condition (1) of Assumption \ref{assu:irrd}, the fiber $h_{n}^{-1}(C)$ of the Hilbert-Chow morphism $h_n:=h_{\beta_n}$ contains rank $n$ semistable vector bundles on $C$. Hence by the deformation-obstruction theory of moduli of sheaves on curves ({\it cf}. \cite[Corollary 4.5.5, p. 114]{HL10}) and the adjunction formula,
\begin{equation}
\label{dimfib}
    \begin{aligned}
     \dim h_{n}^{-1}(C)&\geq n^2(g(C)-1)+1\\
     &=\frac{n^2}{2}(K_S(K_S+K_S))+1=\beta_n^2+1.
    \end{aligned}
\end{equation}
We denote by $U_n$ the open subset of $|\beta_n|$ consisting of smooth curves and by $W_n$ the complement of $U_n$ in $|\beta_n|$.
Combining (\ref{dimks}) and (\ref{dimfib}), it follows that 
\begin{equation}
\label{ineq:irrd}
    \dim h_n^{-1}(W_n)\geq(\chi(\cO_S)+q-2)+(\beta_n^2+1)=\dim h_n^{-1}(U_n)+(q-1),
\end{equation}
where the last identity follows from (\ref{dimsmsupp}).
Now we prove by contradiction. If the assertion were false, then $M_{\beta_n}$ would be irreducible. Since $h_n^{-1}(U_n)$ is nonempty, open and $h_n^{-1}(W_n)$ is a proper closed subset of $M_{\beta_n}$,
$$\dim h_n^{-1}(W_n)<\dim M_{\beta_n}=\dim h_n^{-1}(U_n),$$
which contradicts (\ref{ineq:irrd}).

\end{proof}
We give some explicit examples satisfying Assumption \ref{assu:irrd}.
\begin{exmp}
Let $A$ be a principally polarised Abelian surface and let $\Theta$ be the theta divisor on $A$. Consider the double cover $S \rightarrow A$ branched along a smooth curve  $C_m\in |2m\Theta|$ for $m\in \bZ_{>0}$. Then $S$ is a surface of general type satisfying Assumption \ref{assu:irrd} when $m$ is large.
\end{exmp}
Note that the argument in Proposition \ref{counter} fails if the curve $C$ lies in $|mK_S|$ for $m\ge 2$.   From the general philosophy of asymptotic phenomena in higher rank cases,  we still ask if the following question is true:
\begin{question}
  If $M_\beta$ irreducible when Assumption \ref{assu:irrd} fails and $\beta$ is sufficiently positive?   
\end{question}
\vspace{0.5cm}
\bibliographystyle{alpha}
\bibliography{main}

\begin{thebibliography}{CvGKT20}

\bibitem[ACGH85]{ACGH85}
E.~Arbarello, M.~Cornalba, P.~A. Griffiths, and J.~Harris.
\newblock {\em Geometry of algebraic curves. {V}ol. {I}}, volume 267 of {\em Grundlehren der mathematischen Wissenschaften [Fundamental Principles of Mathematical Sciences]}.
\newblock Springer-Verlag, New York, 1985.

\bibitem[AJ78]{AJ78}
M.~F. Atiyah and J.~D.~S. Jones.
\newblock Topological aspects of {Y}ang-{M}ills theory.
\newblock {\em Comm. Math. Phys.}, 61(2):97--118, 1978.

\bibitem[AK80]{AK80}
Allen~B. Altman and Steven~L. Kleiman.
\newblock Compactifying the {P}icard scheme.
\newblock {\em Adv. in Math.}, 35(1):50--112, 1980.

\bibitem[BBD82]{BBD82}
A.~A. Be\u{\i}linson, J.~Bernstein, and P.~Deligne.
\newblock Faisceaux pervers.
\newblock In {\em Analysis and topology on singular spaces, {I} ({L}uminy, 1981)}, volume 100 of {\em Ast\'{e}risque}, pages 5--171. Soc. Math. France, Paris, 1982.

\bibitem[Bea95]{Beauville95}
Arnaud Beauville.
\newblock Sur la cohomologie de certains espaces de modules de fibr\'{e}s vectoriels.
\newblock In {\em Geometry and analysis ({B}ombay, 1992)}, pages 37--40. Tata Inst. Fund. Res., Bombay, 1995.

\bibitem[Bri77]{Briancon77}
Jo\"{e}l Brian\c{c}on.
\newblock Description de {$H{\rm ilb}\sp{n}C\{x,y\}$}.
\newblock {\em Invent. Math.}, 41(1):45--89, 1977.

\bibitem[CC15]{CC15}
Jinwon Choi and Kiryong Chung.
\newblock The geometry of the moduli space of one-dimensional sheaves.
\newblock {\em Sci. China Math.}, 58(3):487--500, 2015.

\bibitem[CvGKT20]{CGKT20}
Jinwon Choi, Michel van Garrel, Sheldon Katz, and Nobuyoshi Takahashi.
\newblock Local {BPS} invariants: enumerative aspects and wall-crossing.
\newblock {\em Int. Math. Res. Not. IMRN}, 2020(17):5450--5475, 2020.

\bibitem[CW22]{CW22}
Izzet Coskun and Matthew Woolf.
\newblock The stable cohomology of moduli spaces of sheaves on surfaces.
\newblock {\em J. Differential Geom.}, 121(2):291--340, 2022.

\bibitem[dCHM12]{dCHM12}
Mark Andrea~A. de~Cataldo, Tam\'{a}s Hausel, and Luca Migliorini.
\newblock Topology of {H}itchin systems and {H}odge theory of character varieties: the case {$A_1$}.
\newblock {\em Ann. of Math. (2)}, 175(3):1329--1407, 2012.

\bibitem[dCM05]{dCM05}
Mark Andrea~A. de~Cataldo and Luca Migliorini.
\newblock The {H}odge theory of algebraic maps.
\newblock {\em Ann. Sci. \'{E}cole Norm. Sup. (4)}, 38(5):693--750, 2005.

\bibitem[Don90]{Don90}
S.~K. Donaldson.
\newblock Polynomial invariants for smooth four-manifolds.
\newblock {\em Topology}, 29(3):257--315, 1990.

\bibitem[Gie77]{Gieseker77}
D.~Gieseker.
\newblock On the moduli of vector bundles on an algebraic surface.
\newblock {\em Ann. of Math. (2)}, 106(1):45--60, 1977.

\bibitem[GL94]{GL94}
David Gieseker and Jun Li.
\newblock Irreducibility of moduli of rank-{$2$} vector bundles on algebraic surfaces.
\newblock {\em J. Differential Geom.}, 40(1):23--104, 1994.

\bibitem[GL96]{GiesekerLi96}
David Gieseker and Jun Li.
\newblock Moduli of high rank vector bundles over surfaces.
\newblock {\em J. Amer. Math. Soc.}, 9(1):107--151, 1996.

\bibitem[GV98]{GV}
Rajesh Gopakumar and Cumrun Vafa.
\newblock M-theory and topological strings--ii.
\newblock {\em arXiv preprint hep-th/9812127}, 1998.

\bibitem[Gö90]{Got90}
Lothar Göttsche.
\newblock The {B}etti numbers of the {H}ilbert scheme of points on a smooth projective surface.
\newblock {\em Math. Ann.}, 286(1-3):193--207, 1990.

\bibitem[Har77]{Hartshorne}
Robin Hartshorne.
\newblock {\em Algebraic geometry}.
\newblock Graduate Texts in Mathematics, No. 52. Springer-Verlag, New York-Heidelberg, 1977.

\bibitem[Hat02]{Hatcher02}
Allen Hatcher.
\newblock {\em Algebraic topology}.
\newblock Cambridge University Press, Cambridge, 2002.

\bibitem[HL10]{HL10}
Daniel Huybrechts and Manfred Lehn.
\newblock {\em The geometry of moduli spaces of sheaves}.
\newblock Cambridge Mathematical Library. Cambridge University Press, Cambridge, second edition, 2010.

\bibitem[HST01]{HST}
Shinobu Hosono, Masa-Hiko Saito, and Atsushi Takahashi.
\newblock Relative lefschetz action and bps state counting.
\newblock {\em International Mathematics Research Notices}, 2001(15):783--816, 2001.

\bibitem[KL12]{KL12}
Young-Hoon Kiem and Jun Li.
\newblock Categorification of donaldson-thomas invariants via perverse sheaves.
\newblock {\em arXiv preprint arXiv:1212.6444}, 2012.

\bibitem[Kle05]{Kleiman05}
Steven~L. Kleiman.
\newblock The {P}icard scheme.
\newblock In {\em Fundamental algebraic geometry}, volume 123 of {\em Math. Surveys Monogr.}, pages 235--321. Amer. Math. Soc., Providence, RI, 2005.

\bibitem[KLMP24]{KLMP24}
Yakov {Kononov}, Woonam {Lim}, Miguel {Moreira}, and Weite {Pi}.
\newblock Cohomology rings of the moduli of one-dimensional sheaves on the projective plane.
\newblock {\em arXiv preprint arXiv:2403.06277}, 2024.

\bibitem[KPS23]{KPS}
Yakov Kononov, Weite Pi, and Junliang Shen.
\newblock Perverse filtrations, {C}hern filtrations, and refined {BPS} invariants for local {$\Bbb{P}^2$}.
\newblock {\em Adv. Math.}, 433:Paper No. 109294, 29, 2023.

\bibitem[Laz04]{Lazars04}
Robert Lazarsfeld.
\newblock {\em Positivity in algebraic geometry. {I}}, volume~48 of {\em Ergebnisse der Mathematik und ihrer Grenzgebiete. 3. Folge. A Series of Modern Surveys in Mathematics [Results in Mathematics and Related Areas. 3rd Series. A Series of Modern Surveys in Mathematics]}.
\newblock Springer-Verlag, Berlin, 2004.
\newblock Classical setting: line bundles and linear series.

\bibitem[Li97]{Li97}
Jun Li.
\newblock The first two {B}etti numbers of the moduli spaces of vector bundles on surfaces.
\newblock {\em Comm. Anal. Geom.}, 5(4):625--684, 1997.

\bibitem[Lua23]{Luan23}
Yuze Luan.
\newblock Irreducible components of {H}ilbert scheme of points on non-reduced curves.
\newblock {\em arXiv preprint arXiv:2210.01170v2}, 2023.

\bibitem[Mar78]{Maruyama78}
Masaki Maruyama.
\newblock Moduli of stable sheaves. {II}.
\newblock {\em J. Math. Kyoto Univ.}, 18(3):557--614, 1978.

\bibitem[Mar02]{markman02}
Eyal Markman.
\newblock Generators of the cohomology ring of moduli spaces of sheaves on symplectic surfaces.
\newblock {\em J. Reine Angew. Math.}, 544:61--82, 2002.

\bibitem[Mar07]{markman07}
Eyal Markman.
\newblock Integral generators for the cohomology ring of moduli spaces of sheaves over {P}oisson surfaces.
\newblock {\em Adv. Math.}, 208(2):622--646, 2007.

\bibitem[Max19]{Maxim19}
Lauren\c{t}iu~G. Maxim.
\newblock {\em Intersection homology \& perverse sheaves---with applications to singularities}, volume 281 of {\em Graduate Texts in Mathematics}.
\newblock Springer, Cham, [2019] \copyright 2019.

\bibitem[MS13]{MS2013}
Luca Migliorini and Vivek Shende.
\newblock A support theorem for {H}ilbert schemes of planar curves.
\newblock {\em J. Eur. Math. Soc. (JEMS)}, 15(6):2353--2367, 2013.

\bibitem[MS23]{MS23}
Davesh Maulik and Junliang Shen.
\newblock Cohomological {$\chi$}-independence for moduli of one-dimensional sheaves and moduli of {H}iggs bundles.
\newblock {\em Geom. Topol.}, 27(4):1539--1586, 2023.

\bibitem[MSY23]{MSY23}
Davesh {Maulik}, Junliang {Shen}, and Qizheng {Yin}.
\newblock {Perverse filtrations and Fourier transforms}.
\newblock {\em arXiv preprint arXiv:2308.13160}, 2023.

\bibitem[MT18]{MT18}
Davesh Maulik and Yukinobu Toda.
\newblock Gopakumar-{V}afa invariants via vanishing cycles.
\newblock {\em Invent. Math.}, 213(3):1017--1097, 2018.

\bibitem[NS97]{NS97}
D.~S. Nagaraj and C.~S. Seshadri.
\newblock Degenerations of the moduli spaces of vector bundles on curves. {I}.
\newblock {\em Proc. Indian Acad. Sci. Math. Sci.}, 107(2):101--137, 1997.

\bibitem[O'G96]{O'Grady96}
Kieran~G. O'Grady.
\newblock Moduli of vector bundles on projective surfaces: some basic results.
\newblock {\em Invent. Math.}, 123(1):141--207, 1996.

\bibitem[PS08]{PS08}
Chris A.~M. Peters and Joseph H.~M. Steenbrink.
\newblock {\em Mixed {H}odge structures}, volume~52 of {\em Ergebnisse der Mathematik und ihrer Grenzgebiete. 3. Folge. A Series of Modern Surveys in Mathematics [Results in Mathematics and Related Areas. 3rd Series. A Series of Modern Surveys in Mathematics]}.
\newblock Springer-Verlag, Berlin, 2008.

\bibitem[PS23]{PS23}
Weite Pi and Junliang Shen.
\newblock Generators for the cohomology ring of the moduli of 1-dimensional sheaves on {$\Bbb P^2$}.
\newblock {\em Algebr. Geom.}, 10(4):504--520, 2023.

\bibitem[PS24]{PS24}
Weite {Pi} and Junliang {Shen}.
\newblock {On the $P=C$ conjecture and refined BPS invariants for local $\mathbb{P}^2$}.
\newblock {\em arXiv e-prints}, page arXiv:2406.10004, June 2024.

\bibitem[PT10]{PT10}
R.~Pandharipande and R.~P. Thomas.
\newblock Stable pairs and {BPS} invariants.
\newblock {\em J. Amer. Math. Soc.}, 23(1):267--297, 2010.

\bibitem[Ren18]{Rennemo18}
J\o rgen~Vold Rennemo.
\newblock Homology of {H}ilbert schemes of points on a locally planar curve.
\newblock {\em J. Eur. Math. Soc. (JEMS)}, 20(7):1629--1654, 2018.

\bibitem[Sac19]{Sac19}
Giulia Sacc\`a.
\newblock Relative compactified {J}acobians of linear systems on {E}nriques surfaces.
\newblock {\em Trans. Amer. Math. Soc.}, 371(11):7791--7843, 2019.

\bibitem[She12]{She12}
Vivek Shende.
\newblock Hilbert schemes of points on a locally planar curve and the {S}everi strata of its versal deformation.
\newblock {\em Compos. Math.}, 148(2):531--547, 2012.

\bibitem[Sim94]{Simpson94}
Carlos~T. Simpson.
\newblock Moduli of representations of the fundamental group of a smooth projective variety. {I}.
\newblock {\em Inst. Hautes \'{E}tudes Sci. Publ. Math.}, (79):47--129, 1994.

\bibitem[Tau84]{Taubes84}
Clifford~Henry Taubes.
\newblock Path-connected {Y}ang-{M}ills moduli spaces.
\newblock {\em J. Differential Geom.}, 19(2):337--392, 1984.

\bibitem[Yua23a]{Yuan23a}
Yao Yuan.
\newblock On the perverse filtration of the moduli spaces of 1-dimensional sheaves on $\mathbb{P}^2$ and {$P= C$} conjecture.
\newblock {\em arXiv preprint arXiv:2312.17035}, 2023.

\bibitem[Yua23b]{Yuan23}
Yao Yuan.
\newblock Sheaves on non-reduced curves in a projective surface.
\newblock {\em Sci. China Math.}, 66(2):237--250, 2023.

\end{thebibliography}
\end{document}